\documentclass[a4paper,10pt]{article}
\usepackage{amssymb}
\newtheorem{Th}{Theorem}
\newtheorem{Prop}{Proposition}

\newtheorem{Lm}{Lemma}
\newtheorem{Lma}{Lemma}[section]

\newcommand{\be}{\begin{equation}}
\newcommand{\ee}{\end{equation}}
\newcommand{\R}{\mathbb{R}}
\newcommand{\N}{\mathbb{N}}
\newcommand{\C}{\mathbb{C}}

\newcommand\res{\mathop{\hbox{\vrule height 7pt width .5pt depth 0pt
\vrule height .5pt width 6pt depth 0pt}}\nolimits}

\newcommand{\reset}{\setcounter{equation}{0}\setcounter{Th}{0}\setcounter{Prop}{0}\setcounter{Co}{0}
\setcounter{Lm}{0}\setcounter{Rm}{0}}
\def\La{\Lambda}
\def\La{\Lambda}
\def\ti{\tilde}
\def\lf{\left}
\def\rg{\right}

\def\al{\alpha}
\def\la{\lambda}

\def\ep{\varepsilon}
\def\ds{\displaystyle}
\def\ov{\overline}
\def\Om{\Omega}
\def\om{\omega}
\def\p{\partial}
\def\bn{\vec{n}}

\def\bbe{\vec{e}}
\def\bbf{\vec{f}}

\def\bL{\vec{L}}
\def\bR{\vec{R}}

\def\res{\mathop{\hbox{\vrule height 7pt width .5pt 
depth 0pt\vrule height .5pt width 6pt depth 0pt}}\nolimits}
\begin{document}
\title{Energy Quantization for Willmore Surfaces\\
 and Applications}
\author{Yann Bernard\footnote{Mathematisches Institut, Albert-Ludwigs-Universit\"at, 79004 Freiburg, Germany. Supported by the DFG SFB 71 fund (project B3).}\,\, and Tristan Rivi\`ere\footnote{Department of Mathematics, ETH Zentrum,
CH-8093 Z\"urich, Switzerland.}}
\date{ }
\maketitle

{\bf Abstract :} {\it We prove a bubble-neck decomposition together with an energy quantization result for sequences of Willmore surfaces into ${\R}^m$ with uniformly bounded energy and non-degenerating conformal type. We deduce the strong compactness of Willmore closed surfaces of a given genus modulo the M\"obius group action, below some energy threshold.}

\medskip

\noindent{\bf Math. Class.} 30C70, 58E15, 58E30, 49Q10, 53A30, 35R01, 35J35, 35J48, 35J50.

\section{Introduction}

Let $\vec{\Phi}$ be an immersion from a closed abstract two-dimensional manifold $\Sigma$ into ${\R}^{m\ge3}$. We denote by $g:=\vec{\Phi}^\ast g_{{\R}^m}$ the pull back by $\vec{\Phi}$ of the flat  canonical metric $g_{{\R}^m}$ of ${\R}^m$, also called the {\it first fundamental form of $\vec{\Phi}$}, and we let $dvol_g$ be its associated {\it volume form}. The {\it Gauss map} of the immersion $\vec{\Phi}$ is the map taking values in the Grassmannian
of oriented $m-2$-planes in ${\R}^m$ given by
\[
\vec{n}_{\vec{\Phi}}:=\star\frac{\p_{x_1}\vec{\Phi}\wedge\p_{x_2}\vec{\Phi}}{|\p_{x_1}\vec{\Phi}\wedge\p_{x_2}\vec{\Phi}|}\quad,
\] 
where $\star$ is the usual Hodge star operator in the Euclidean metric. \\
Denoting by $\pi_{\vec{n}_{\vec{\Phi}}}$ the orthonormal projection of vectors in ${\R}^m$ onto the $m-2$-plane given by $\vec{n}_{\vec{\Phi}}$, the {\it second fundamental form} may be expressed as\footnote{In order to define $d^2\vec{\Phi}(X,Y)$ one has to extend locally around $T_p\Sigma$ the vector-fields $X$ and $Y$. It is not difficult to
check that $\pi_{\vec{n}_{\vec{\Phi}}}d^2\vec{\Phi}(X,Y)$ is independent of this extension.}

\[
\forall\:\: X,Y\in T_p\Sigma\quad\quad\quad \vec{\mathbb I}_p(X,Y):=\pi_{\vec{n}_{\vec{\Phi}}}d^2\vec{\Phi}(X,Y)\quad.
\]The {\it mean curvature vector} of the immersion at $p$ is 
\[
\vec{H}_{\vec{\Phi}}:=\frac{1}{2}\,tr_g(\vec{\mathbb I})=\frac{1}{2}\,\lf[\vec{\mathbb I}(\ep_1,\ep_1)+\vec{\mathbb I}(\ep_2,\ep_2)\rg]\quad,
\]
where $(\ep_1,\ep_2)$ is an orthonormal basis of $T_p\Sigma$ for the metric ${g_{\vec{\Phi}}}$. 

\medskip

In the present paper, we study the Lagrangian given by the $L^2$-norm of the second fundamental form:
\[
E(\vec{\Phi}):=\int_{\Sigma} |\vec{\mathbb I}|^2_{g_{\vec{\Phi}}}\ dvol_{g_{\vec{\Phi}}}\quad.
\]
An elementary computation gives 
\[
E(\vec{\Phi}):=\int_{\Sigma} |\vec{\mathbb I}|^2_{g_{\vec{\Phi}}}\ dvol_{g_{\vec{\Phi}}}=\int_{\Sigma} |d\vec{n}_{\vec{\Phi}}|^2_{g_{\vec{\Phi}}}\ dvol_{g_{\vec{\Phi}}}\quad.
\]
The energy $E$ may accordingly be seen as the {\it Dirichlet Energy} of the Gauss map $\vec{n}_{\vec{\Phi}}$ with respect to the induced
metric $g_{\vec{\Phi}}$.
The Gauss Bonnet theorem implies that 
\be
\label{0.1}
E(\vec{\Phi}):=\int_{\Sigma} |\vec{\mathbb I}|^2_{g_{\vec{\Phi}}}\ dvol_{g_{\vec{\Phi}}}=4\int_{\Sigma}|\vec{H}_{\vec{\Phi}}|^2\ dvol_{g_{\vec{\Phi}}}-2\int_{\Sigma} K_{{\vec{\Phi}}}\ dvol_{g_{\vec{\Phi}}}=4\ \int_{\Sigma}|\vec{H}_{\vec{\Phi}}|^2\ dvol_g-4\pi\chi(\Sigma)\quad,
\ee
where $ K_{{\vec{\Phi}}}$ is the Gauss curvature of the immersion, and $\chi(\Sigma)$ is the {\it Euler characteristic} of the surface $\Sigma$.
The energy
\[
W(\vec{\Phi}):=\int_\Sigma|\vec{H}_{\vec{\Phi}}|^2\ dvol_{g_{\vec{\Phi}}}\quad,
\]
is called {\it Willmore energy}.

\medskip

Critical points of the {\it Willmore energy}, comprising for example minimal surfaces\footnote{minimal  surfaces satisfy $\vec{H}=0$ and are hence absolute minimizers of $W$.}, are called {\it Willmore surfaces.} Although already known in the XIXth century in the context of the elasticity theory of plates, it was first considered in conformal geometry by Blaschke in \cite{Bla3} who sought to merge
the {\it theory of minimal surfaces} and the {\it conformal invariance property}. This Lagrangian has indeed both desired features : its critical points contain minimal surfaces, and it is conformal invariant, owing to the following pointwise identity which holds for an arbitrary immersion $\vec{\Phi}$ of $\Sigma$ into ${\R}^m$ and at every point of $\Sigma$\,:
\[
\forall\:\ \Xi \quad\mbox{ conformal diffeo. of }{\R}^m\cup\{\infty\} \quad \lf[|\vec{H}_{\Xi\circ\vec{\Phi}}|^2-K_{\Xi\circ{\vec{\Phi}}}\rg] dvol_{g_{\Xi\circ\vec{\Phi}}}=\lf[|\vec{H}_{\vec{\Phi}}|^2-K_{{\vec{\Phi}}}\rg] dvol_{g_{\vec{\Phi}}}\quad.
\]
Using again Gauss Bonnet theorem, the latter implies the conformal invariance of $W$\,:
\[
\forall\:\ \Xi \quad\mbox{ conformal diffeo. from }{\R}^m\cup\{\infty\} \quad\quad W(\Xi\circ\vec{\Phi})=W(\vec{\Phi})\quad.
\]
This conformal invariance implies that the image of a {\it Willmore immersion} by a conformal transformation of ${\R}^m$ is still a Willmore immersion. Starting for example from a minimal surface, one may then generate many new Willmore surfaces, simply by applying conformal transformations (naturally, these surfaces need no longer be minimal). In his time, Blaschke used the term {\it conformal minimal} for the critical points of $W$, seeking to insist on this idea of merging together the theory of minimal surface with conformal invariance.

\medskip

An important task in the analysis of Willmore surfaces is to understand the closure of the space of {\it Willmore immersions} under a certain level of energy. Because of the non-compactness of the conformal group of transformation of ${\R}^m$, one cannot expect that the space of Wilmore immersions in closed in the strong $C^l$-topology. However,
locally, in isothermic coordinates\footnote{Analogously to other gauge-invariant problems, such as in Yang-Mills theory, isothermic coordinates or conformal parametrizations provide the optimal symmetry breaking method in the search of pertinent estimates. A detailed discussion on this topic is available in \cite{Ri1}.}, under some universal energy threshold, if the conformal factor of the induced metric $g_{\vec{\Phi}}$ is controlled in $L^\infty$, then 
the immersion is uniformly bounded in any $C^l$ norm. More precisely there holds the following $\ep$-regularity result.

\begin{Th}
\label{th-I.0} \cite{Ri2}
There exists $\ep(m)>0$ such that, for any Willmore conformal immersion $\vec{\Phi}$ from $B_1(0)$ into ${\R}^m$ satisfying
\[
\int_{B_1(0)} |\nabla \vec{n}_{\vec{\Phi}}|^2\, dx<\ep(m)
\]
then for any $l\in {\N}^\ast$ we have
\[
 \| e^{-\la}\ \nabla^l\vec{\Phi}\|_{L^\infty(B_{1/2})}\le C_l\  \lf[\int_{B_1(0)} |\nabla \vec{n}_{\vec{\Phi}}|^2\, dx\, + \,1\rg]^{1/2}\quad,
\]
where $C_l$ only depends on $l$, while $\la$ denotes the conformal parameter of $\vec{\Phi}$. Namely, $\la=\|\log|\p_{x_1}\vec{\Phi}|\|_{L^\infty(B_1)}=\|\log|\p_{x_2}\vec{\Phi}|\|_{L^\infty(B_1)}$.\hfill $\Box$
\end{Th}

This theorem leads to the {\it concentration of compactness} ``dialectic" developed by Sacks and Uhlenbeck. In a conformal parametrization, assuming that the conformal factor is $L^\infty$-controlled in some subdomain of $\Sigma$, then a sequence of Willmore immersions can fail to convergence strongly in $C^l$ only at isolated points ; namely, at those points where the $W^{1,2}$-norm of the Gaus map concentrates.
Assuming their induced metric  generates a sequence of conformal classes which remains within a compact subdomain of the moduli space of $\Sigma$, the control of the conformal factor of a sequence of conformal immersions with uniformly bounded Willmore energy is also guaranteed, except at those isolated points of $\Sigma$ where the $W^{1,2}$-norm of the Gaus map concentrates.
This fact is established in \cite{Ri3} (see also proposition~\ref{pr-III.1} below), and it utlimately follows from the works of Toro \cite{To}, of M\"uller and Sverak \cite{MS}, and from the work of H\'elein \cite{Hel} on immersions with totally bounded curvature.

\medskip

To fully understand the loss of strong compactness of a sequence of Willmore surfaces, the difficulty is of course to perform the necessary analysis at the ``blow up points". This is the aim of the present paper.
Under the assumption that the sequence of conformal classes associated to the sequence of Willmore surfaces remains within a compact subdomain of the Moduli space, we show that, modulo extraction of a subsequence, the parametrizing abstract surface splits into three distinct regions : the main region where strong convergence holds, the concentrating  parametrization of non-trivial Willmore spheres, and finally the bubble and neck regions connecting the two previous ones and in which we show the energy vanishes.

\medskip

We now state our main result.

\begin{Th}
\label{th-I.1}
Let $\vec{\Phi}_k$ be a sequence of Willmore immersions of a closed surface $\Sigma$. Assume that
\[
\limsup_{k\rightarrow +\infty} W(\vec{\Phi}_k)<+\infty\quad,
\]
that the conformal class of $\vec{\Phi}_k^\ast g_{{\R}^m}$ remains within a compact subdomain of the moduli space of $\Sigma$. 
Then, modulo extraction of a subsequence, the following energy identity holds
\be
\label{0.2}
\lim_{k\rightarrow +\infty}W(\vec{\Phi}_k)=W(\vec{\xi}_\infty)+\sum_{s=1}^p W(\vec{\eta}_s)+\sum_{t=1}^q\lf[W(\vec{\zeta}_t)-4\pi\,\theta_t\rg]\quad,
\ee
where $\vec{\xi}_\infty$ is a possibly branched smooth immersion of $\Sigma$. The maps $\vec{\eta}_s$ and $\vec{\zeta}_t$ are smooth, possibly branched, immersions
of $S^2$ ; and $\theta_t$ is the integer density of the current $(\vec{\zeta}_t)_\ast[S^2]$ at some point $p_t\in \vec{\zeta}_t(S^2)$, namely
\[
\theta_t:=\lim_{r\rightarrow 0} \frac{{\mathcal H}^2\lf(B^m_r(p_t)\cap\vec{\zeta}_t(S^2)\rg)}{\pi\,r^2}\quad.
\] 
\end{Th}
The second part of our main result describes how the maps $\vec{\xi}_\infty$, $\vec{\eta}^s$ and $\vec{\zeta}^t$ are obtained from the original sequence $\vec{\Phi}_k$. 

\begin{Th}
\label{th-I.2}
With the same notation as in theorem~\ref{th-I.1}, the immersion $\vec{\xi}_\infty$ is obtained as follows. There exist a sequence $f_k$ of diffeomorphisms of $\Sigma$\,, a sequence $\vec{\Xi}_k$
of M\"oebius transformations of ${\R}^m$\,, and finitely many points $\{a^1,\ldots, a^n\}$ such that $\Xi_k\circ\vec{\xi}_k\circ f_k$ is conformal and
\be
\label{0.3}
\vec{\xi}_k:=\Xi_k\circ\vec{\xi}_k\circ f_k\longrightarrow\vec{\xi}_\infty\quad\quad\mbox{ in }C^l_{loc}(\Sigma\setminus\{a^1,\ldots, a^n\})\quad\forall\:\ l\in {\N}\quad.
\ee
Furthermore, there holds
\be
\label{0.4}
\lim_{k\rightarrow +\infty}W(\vec{\Phi}_k)=W(\vec{\xi}_\infty)\quad\quad\Longleftrightarrow\quad\quad \vec{\xi}_k\longrightarrow\vec{\xi}_\infty\quad\quad\mbox{ in }C^l(\Sigma)\quad\forall\:\ l\in {\N}\quad.
\ee
Finally there exists a sequence $h_k$ of constant scalar curvature metric conformally equivalent to $\vec{\xi}_k=\Xi_k\circ\vec{\xi}_k\circ f_k$ and strongly converging
in $C^l(\Sigma)$, such that for any $s\in\{1,\ldots, p\}$ (resp. for any $t\in\{1,\ldots, q\}$) there exist a sequence of points $x^s_k\in \Sigma$ converging to one of the $a^i$\,, a sequence of radius $\rho^s_k$ converging
to zero, and a sequence of M\"obius transformations $\Xi_k^s$ (resp. $\Xi_k^t$) for which (in converging $h_k$ conformal coordinates $\varphi_k$ around $a^i$), one has
\[
\Xi_k^s\circ\vec{\xi}_k\circ\varphi_k( \rho_k^s\ y+\phi^{-1}_k(x_k^s))\longrightarrow \vec{\eta}_s\circ \pi^{-1}(y)\quad\quad\quad\mbox{ in }\quad C^l_{loc}({\C}\setminus 
\{a_{s}^1,\ldots, a_s^{n_s}\})
\]
and respectively
\[
\Xi_k^t\circ\vec{\xi}_k\circ\varphi_k( \rho_k^t\ y+\phi^{-1}_k(x_k^t))\longrightarrow \vec{\zeta}_t\circ \pi^{-1}(y)\quad\quad\quad\mbox{ in }\quad C^l_{loc}({\C}\setminus 
\{a_{t}^1,\ldots, a_t^{n_t}\})
\]
for any $l\ \in {\N}$ ; where $\pi$ denotes the stereographical projection from $S^2$ into ${\C}$ and $\{a_{s}^1,\ldots, a_s^{n_s}\}$ and $\{a_{t}^1,\ldots, a_t^{n_t}\}$ are finite sets of points in the complex plane.\\
The maps $\Xi_k^s$ are compositions of a dilation and isometries, while each map $\Xi_k^t$ is a composition of one inversion, one dilation, and one isometry. 
\hfill $\Box$
\end{Th}

The results given in theorem~\ref{th-I.1} and theorem~\ref{th-I.2} are to be viewed in the context of other bubble-neck decomposition and energy quantization situations previously studied. In particular for harmonic maps and other conformally invariant problems as those considered in \cite{SaU}, \cite{Jo}, \cite{St}, \cite{DiT}, \cite{Pa}, \cite{Ri5}, and in\cite{Ri6}. However, these problems are all of second-order elliptic or parabolic types. The novelty of the present work is an energy quantization result and bubble-neck decomposition for a fourth-order problem\footnote{Indeed, in a conformal parametrization $\vec{\xi}$, the Willmore functional may be recast as
\[
W(\vec{\xi})=\frac{1}{4}\int_{\Sigma}|\Delta_g\vec{\xi}|^2\ dvol_g\quad,
\]
thereby giving rise to a fourth-order problem}. Moreover, the analysis which we present here requires to handle one additional delicate difficulty, namely the need to first control the conformal parameter in ``conformally degenerating'' neck regions (annuli), prior to obtaining from PDE techniques the vanishing of the energy in these regions.

\medskip

The method we develop to prove theorem~\ref{th-I.1} and theorem~\ref{th-I.2}  is partly inspired by a technique
introduced by the second author in collaboration with F.H. Lin, in the context of the Ginzburg-Landau equations \cite{LiRi1} and in the context of harmonic maps  \cite{LiRi2}. This techniques relies on the properties of the Lorentz interpolation spaces. More precisely, in the present paper, we obtain estimates of the trace of the second fundamental form in the neck regions, respectively in the weak Marcinkiewicz space $L^{2,\infty}$ space, and in its predual, the Lorentz space $L^{2,1}$. These estimates constitute the main achievement of our proof. 
Prior to deriving these estimates, it is first necessary to carefully study the conformal factor of the conformal
Willmore immersions in the neck region. This is done through extracting a good Coulomb gauge. The existence of such a gauge is one of the novelties of the present work, for it requires only the gradient of the Gauss map lie in $L^2$, and that it have a small $L^{2,\infty}$-norm.
Such a control on the $L^{2,\infty}$-norm of the second fundamental form (i.e. of the gradient of the Gauss map) ultimately  stems from the $ep$-regularity established in \cite{Ri2}. The appropriate control in $L^{2,1}$ is then deduced using the standard methods of integrability by compensation and applying them to the conservation laws satisfied by Willmore immersions, and originally derived in \cite{Ri2}.

\medskip

It is legitimate to ask what happens if one removes the assumption on the control of the conformal structure. In this situation, the first important observation is provided by the examples constructed in \cite{BuPa}, where inverted catenoids arise as limiting bubbles, which are thus not smooth Willmore immersion of $S^2$. The second important observation is that harmonic mapping from a degenerating Riemann surface into a manifold are known in general not to
satisfy the energy quantization property in the thin collar regions (cf. example constructed in \cite{Pa}, and the careful systematic study given in \cite{Zhu}). Questions associated with the more complex situation in which the sequence of conformal classes no longer remains within a subdomain of the Moduli space are considered in a forthcoming work \cite{BR3}.

\medskip

To prevent the conformal class from degenerating, a sufficient energy condition has been given (in the cases $m=3$ and $m=4$) in \cite{KS1}, and for surfaces in any codimension in \cite{Ri4} (and independently in \cite{KuLi}).  Following the convention in \cite{Si}, we introduce
 \[
 \beta_g^m:=\inf\lf\{W(\vec{\Phi})\quad;\quad \vec{\Phi}\mbox{ is an immersion of the genus }g\mbox{ closed surface }\rg\}\quad. 
 \]
 and 
 \[
 \om_g^m:=\min\lf\{4\pi+\sum_{i=1}^p(\beta_{g_i}^m-4\pi)\quad\quad;\quad g=g_1+\ldots+ g_p\quad, \quad 1\le g_i< g  \rg\}\quad.
 \]
In \cite{BK}, it is proved that for $g\ge 2$ there holds
\be
\label{aa1}
\beta_g^m<\om_g^m\quad.
\ee 
 \begin{Th}
\label{pr-III.4} \cite{Ri4} \cite{KuLi}
Let $(\Sigma,c_k)$ be a sequence of closed Riemann surfaces of genus $g$, with degenerating conformal class $[c_k]$ diverging to the boundary of the Moduli space of $\Sigma$. Let $\vec{\Phi}_k$ be a sequence of conformal immersions
from $(\Sigma,c_k)$ into ${\R}^m$. Then 
\be
\label{III.u0}
\liminf_{k\rightarrow +\infty}\int_{\Sigma}|\vec{H}_{\vec{\Phi}_k}|^2\ dvol_{\vec{\Phi}_k^\ast g_{{\R}^m}}\ge \min\{8\pi,\om^m_g\}\quad.
\ee
\hfill $\Box$
\end{Th}
A consequence of our analysis below is the following result.
\begin{Th}
\label{th-I.3}
Let $\Sigma$ be an arbitrary closed two-dimensional manifold. Modulo the action of the M\"obius group of ${\R}^m$, the space of Willmore immersions into ${\R}^m$, for $m=3$ and $m=4$, satisfying
\[
W(\Phi)<\min\{8\pi,\om^m_g\}-\delta
\]
is strongly compact in the $C^l$ topology for any $l\in {\N}$ and any $\delta>0$.\hfill $\Box$
\end{Th}
This result was obtained for $\Sigma=T^2$ and $m=3$ in \cite{KS} and for $\Sigma=T^2$ and $m=4$ in \cite{Ri2}.

\medskip\medskip

The paper is organized as follows. In section II, we recall the formalism introduce in \cite{Ri3} which allow us to appropriately ``renormalize" a sequence of weak immersions uniformly bounded Willmore energy. 
In section III, we outline a generic ``energy tracking" procedure enabling us to detect bubbles and neck regions. Section IV, also general in nature, explains how to construct a Coulomb moving frame under the hypothesis that the Gauss map has small weak-$L^2$-norm. In Section V, we develop uniform Harnack estimates giving control of the conformal factor of weak immersions in the neck regions. Willmore immersions first appear in section VI, where we show how to control the $L^{2,1}$-norm of the mean curvature vector in the neck regions. In Section VII, we derive an energy quantization result for a sequence of Gauss maps corresponding to a sequence of Willmore immersions. Finally, in section VIII, all of the aforementioned estimates are seamed together, and the main theorems~\ref{th-I.1} and \ref{th-I.2} are proved.

\medskip

\noindent{\bf Acknowledgements}. This work was partly elaborated while the first author was visiting the Forschungsinstitut f\"ur Mathematik at the ETH Z\"urich.
He would like to thank the institute for it's hospitality and the excellent working conditions.

\section{Preliminaries}
\reset
In \cite{Ri3}, the second author developed a suitable framework in which one can perform the calculus of variation of the Willmore Lagrangian of a surface. In this framework, a particularly useful space is that of Lipschitz immersions with $L^2-$bounded second fundamental forms.\\
Let $g_0$ be a smooth ``reference" metric on $\Sigma$. One defines the Sobolev spaces $W^{k,p}(\Sigma,{\R}^m)$ of measurable maps from $\Sigma$ into 
${\R}^m$ as
\[
W^{k,p}(\Sigma,{\R}^m):=\lf\{f\ \mbox{ meas. } {\Sigma}\rightarrow {\R}^m\mbox{ s.t. }\sum_{l=0}^k\int_{\Sigma}|\nabla^l f|_{g_0}^p\ dvol_{g_0}<+\infty\rg\}\quad.
\]
As $\Sigma$ is assumed to be compact, this space is independent of the chosen reference metric $g_0$.

\medskip

It is important to have a weak first fundamental form: we need  $\vec{\Phi}^\ast g_{{\R}^m}$ to define an $L^\infty$ metric with a bounded inverse.
This is the case if we assume that $\vec{\Phi}$ lies in $W^{1,\infty}(\Sigma)$ and if $d\vec{\Phi}$ has maximal rank 2 at every point, with some uniform quantitative control
of ``how far'' $d\vec{\Phi}$ is from being degenerate. Namely, there exists $c_0>0$ such that
\be
\label{I.1}
|d\vec{\Phi}\wedge d\vec{\Phi}|_{g_0}\ge c_0>0 \quad,
\ee
where $d\vec{\Phi}\wedge d\vec{\Phi}$ is a 2-form on $\Sigma$ with values in the space of 2-vectors of ${\R}^m$, and given in local coordinates
by $2\,\p_x\vec{\Phi}\wedge\p_y\vec{\Phi}\ dx\wedge dy$.         
Note that the condition (\ref{I.1}) is independent of the choice of the metric $g_0$ .
To a Lipschitz immersion satisfying (\ref{I.1}) we associate its Gauss map
\[
\vec{n}_{\vec{\Phi}}:=\star\frac{\p_{x_1}\vec{\Phi}\wedge\p_{x_2}\vec{\Phi}}{|\p_{x_1}\vec{\Phi}\wedge\p_{x_2}\vec{\Phi}|}\quad.
\] 
It belongs to $L^\infty(\Sigma)$ and takes values in the Grassmannian
of oriented $(m-2)$-planes in ${\R}^m$.\\

\noindent
We next introduce the space ${\mathcal E}_{\Sigma}$ of Lipschitz immersions\footnote{henceforth called {\it weak immersions.} } of $\Sigma$ 
with square-integrable second fundamental form:
\[
\mathcal{E}_\Sigma:=\lf\{
\begin{array}{c}
\ds\vec{\Phi}\in W^{1,\infty}(\Sigma,{\R}^m)\quad\mbox{ s.t. } \vec{\Phi} \mbox{ satisfies (\ref{I.1}) for some }c_0\\[5mm]
\ds\mbox{ and }\quad\int_{\Sigma}|d\vec{n}|_g^2\ dvol_g<+\infty
\end{array}
\rg\}\quad .
\]
As before, $g:=\vec{\Phi}^\ast g_{{\R}^m}$ is the pull-back by $\vec{\Phi}$ of the flat canonical metric $g_{{\R}^m}$ of ${\R}^m$, and $dvol_g$ is its corresponding volume form.

\medskip
%
%
%
%
%
%
%
%
%
%
%
%
%
%
%
%
%
\medskip  

As in the case of the Yang-Mills functional, one of the main difficulties one encounters in the variational study of the Willmore functional is the size of the corresponding gauge group (i.e. the group of transformations through which the energy is invariant). For the energy $W$, this gauge group is very large: it contains all diffeomorphisms of $\Sigma$. To ``break'' the gauge invariance, one of the first techniques considered in \cite{Ri3} consists in working with conformal immersions. Indeed, owing to the work of M\"uller and Sverak \cite{MS} and to the work of  H\'elein \cite{Hel}, one gains an $L^\infty$-control of the conformal parameter under a small energy $E$ assumption. It is then possible to show that the pull-back metric of any element in ${\mathcal E}_\Sigma$ defines a smooth conformal structure on $\Sigma$ (cf. details in \cite{Ri1}).
Thus, corresponding to any element $\vec{\Phi}\in {\mathcal E}_\Sigma$, there exists a Lipschitz diffeomorphism $f$ of $\Sigma$ such that $\vec{\Phi}\circ f$ is conformal.
Having found a conformal structure $c$ of $\Sigma$, one defines
\[
\mathcal{E}^c_\Sigma:=\lf\{
\ds\vec{\Phi}\in {\mathcal E}_\Sigma\quad\mbox{ s.t. }\vec{\Phi}\mbox{ is conformal w.r.t. }c\ \rg\}\quad .
\]
Naturally, in considering a sequence of elements in ${\mathcal E}_\Sigma$, the corresponding sequence of conformal classes it yields is not a priori controlled, and it might diverge in the moduli space of $\Sigma$, even when the Willmore energy of the immersions is uniformly bounded. 
This possible degeneracy of the conformal
class is another source of trouble when trying to control a sequence of immersions with uniformly bounded Willmore energy. In the present work, we will assume that the aforementioned sequence of conformal classes does \textbf{not} degenerate. Sufficient conditions to avoid degeneracy of the conformal structure were given in \cite{Ri4}, and independently in \cite{KuLi}.

\medskip

Once the global gauge choice is made, the second difficulty one encounters is the invariance of the Willmore functional under conformal transformations of ${\R}^m$. Namely, any  conformal diffeomorphism $\Xi$ of ${\R}^m$ (i.e. an element of the M\"obius group of ${\R}^m$) leaves both ${\mathcal E}^c_\Sigma$ and $W$ (through composition). Since the M\"obius group is non-compact, this invariance generates serious troubles in the study of a sequence of immersions with uniformly bounded Willmore energy. 
Indeed, such a sequence may a-priori collapse to a point in the limit! To overcome this difficulty, a method for extracting a suitable $\Xi$ was devised in Lemma A.4 of \cite{Ri3}. Combined with the result of Lemma III.1 in \cite{Ri3}, one obtains an ``almost'' weak closure result for sequences in ${\mathcal E}_\Sigma^c$ modulo composition with suitable elements of the M\"obius group. More precisely, the limit immersion is not necessarily a weak immersion of the whole $\Sigma$, but only a weak immersion away from possibly finitely many points. We are thus led to consider the space of weak branched immersions,
\[
\mathcal{F}_\Sigma:=\lf\{
\begin{array}{c}
\ds\vec{\Phi}\in W^{1,\infty}(\Sigma,{\R}^m)\quad\mbox{ s.t. } \exists\{a_1,\ldots, a_N\}\in \Sigma\quad\mbox{with}\\[5mm]
\ds\forall\: K\mbox{ compact subset of }\Sigma\setminus\{a_1,\ldots, a_N\}\quad\\[5mm]
\ds\quad\exists\, c_K>0\quad\mbox{ s.t. }\forall\: x\in K\quad|d\vec{\Phi}\wedge d\vec{\Phi}|_{g_0}(x)\ge c_K\\[5mm]
\ds\quad\quad\mbox{ and }\quad\int_{\Sigma}|d\vec{n}|_g^2\ dvol_g<+\infty
\end{array}
\rg\}\quad .
\]
Finally, we introduce the space of elements in ${\mathcal F}_\Sigma$ which are conformal with respect to some fixed conformal structure $c$,
\[
\mathcal{F}^c_\Sigma:=\lf\{\vec{\Phi}\in {\mathcal F}_\Sigma\quad\mbox{ s. t. }\quad\vec{\Phi}\mbox{ is conformal w.r.t. }c\ \rg\}\quad.
\]
The above discussion is summarized in the following result, which will be a central ingredient in the rest of this article. 
\begin{Lm}
\label{lm-I.1} \cite{Ri3}
Let $\Sigma$ be a closed two-dimensional manifold. Let $\vec{\Phi}_k$ be a sequence of elements in ${\mathcal E}_\Sigma$ such that $W(\Phi_k)$ is uniformly bounded.
Assume that the conformal class of the conformal structure $c_k$ (i.e. complex structure of $\Sigma$) defined by $\vec{\Phi}_k$ remains in a compact subspace of the Moduli space of $\Sigma$. Then, modulo
extraction of a subsequence, the sequence $c_k$ converges to a smooth limiting complex structure $c_\infty$ ; and there
exist a sequence of Lipschitz diffeomorphisms $f_k$ of $\Sigma$ and a sequence of smooth conformal structures $c_k$ of $\Sigma$ such that $\vec{\Phi}_k\circ f_k$ is conformal from $(\Sigma,c_k)$ into ${\R}^m$. Moreover, there exists a sequence $\Xi_k$ 
of conformal diffeomorphisms of ${\R}^m\cup\{\infty\}$ and at most finitely many points $\{a^1,\ldots, a^n\}$ such that  
\be
\label{II.0}
\limsup_{k\rightarrow +\infty} {\mathcal H}(\Xi_k\circ\vec{\Phi}_k\circ f_k(\Sigma))<+\infty\quad\quad, \quad\quad\Xi_k\circ\vec{\Phi}_k\circ f_k(\Sigma)\subset B_R(0)
\ee
for some $R>0$ independent of $k$, and
\be
\label{II.1}
\vec{\xi}_k:=\Xi_k\circ\vec{\Phi}_k\circ f_k\;\rightharpoonup \;\vec{\xi}_\infty\quad\quad\mbox{ weakly in }\:(W^{2,2}_{loc}\cap W^{1,\infty}_{loc})^\ast(\Sigma\setminus\{a^1,\ldots, a^n\})\quad.
\ee
The convergences are understood with respect to $h_k$, which is the constant scalar curvature metric of unit volume attached to the conformal structure $c_k$. \\
Furthermore, there holds
\be
\label{II.2}
\forall\: K\mbox{ compact subset of }\Sigma\setminus\{a^1,\ldots, a^n\}\quad\limsup_{k\rightarrow+\infty}\|\log|d\vec{\xi}_k|_{h_k}\|_{L^\infty(K)}<+\infty\quad.
\ee
Finally, $\vec{\xi}_\infty$ is a weak immersion of $\Sigma\setminus\{a^1,\ldots, a^n\}$
and conformal from $(\Sigma,c_\infty)$ into ${\R}^m$.\hfill$\Box$
\end{Lm}

\section{The bubble-neck decomposition procedure}

The goal of this section is to establish the following technical proposition. Although its statement might seem at first somewhat overwhelming, its proof involves no particularly challenging difficulties and is based on an iteration argument.  

\begin{Prop}
\label{pr-III.1}{\bf [Bubble-neck decomposition]}
Let $\Sigma$ be a closed two-dimensional manifold. Let $\vec{\xi}_k$ be a sequence of weak Lipschitz immersions with $L^2-$bounded second fundamental
forms (i.e. $\vec{\xi}_k\in{\mathcal E}_\Sigma$) such that
\[
\int_{\Sigma}|d\vec{n}_{\vec{\xi}_k}|^2_{g_k}\ dvol _{g_k}\le \La\quad,
\]
where $g_k:=\vec{\xi}_k^\ast\, g_{{\R}^m}$. 
\smallskip
\noindent
We suppose that $g_k$ is conformally equivalent to a constant scalar curvature metric $h_k$ such that for all $l\in{\N}$
 \[
 h_k\longrightarrow h_\infty \quad\quad\quad\mbox{ in }C^l(\Sigma)\quad,
 \]
 where $h_\infty$ is a constant scalar curvature on $\Sigma$. 

\smallskip
\noindent
Suppose there exist $n$ points $\{a^1,\ldots, a^n\}\subset\Sigma$ and a radius $R>0$ such that the following holds
\begin{itemize}
\item[(i)]
\be
\label{III.4}
\limsup_{k\rightarrow +\infty} {\mathcal H}(\vec{\xi}_k(\Sigma))<+\infty\quad,
\ee
\item[(ii)]
\be
\label{III.5}
\vec{\xi}_k(\Sigma)\subset B_R(0)\quad,
\ee
\item[(iii)]
\be
\label{III.6}
\vec{\xi}_k\;\rightharpoonup\; \vec{\xi}_\infty\quad\quad\mbox{ weakly in }\:(W^{2,2}_{loc}\cap W^{1,\infty}_{loc})^\ast(\Sigma\setminus\{a^1,,\ldots, a^n\})\quad,
\ee
\item[(iv)]
\be
\label{III.7}
\forall\: K\mbox{ compact subset of }\Sigma\setminus\{a^1,\ldots, a^n\}\quad\limsup_{k\rightarrow+\infty}\|\log|d\vec{\xi}_k|_{h_k}\|_{L^\infty(K)}<+\infty\quad.
\ee
\end{itemize}
Then for any $0<\ep_0<8\pi/3$, there exist a subsequence, still denoted $\vec{\xi}_k$\,,\,$n$ integers $\{Q^1,\ldots, Q^n\}$\,,\, $n$ sequences of points $(x_k^{i,j})_{j=1,\ldots, Q^i}\subset\Sigma\,$, and $n$ sequences of radii $(\rho_k^{i,j})_{j=1,\ldots, Q^i}$ satisfying
\be
\label{III.8}
\forall\:i\in{1,\ldots, n}\quad\forall\: j\:\in\{1,\ldots, Q^i\}\quad\quad\quad\lim_{k\rightarrow +\infty}x_k^{i,j}=a^i\quad,
\ee
\be
\label{III.9}
\forall\: i\in{1,\ldots, n}\quad\forall\: j\:\in\{1,\ldots, Q^i\}\quad\quad\quad\lim_{k\rightarrow +\infty}\rho_k^{i,j}=0\quad,
\ee
and
\be
\label{III.10}
\forall\: i\in{1,\ldots, n}\quad\forall\: j\ne j'\in\{1,\ldots, Q^i\}\quad
\lf\{
\begin{array}{l}
\ds\mbox{ either }\lim_{k\rightarrow +\infty}\frac{\rho^{i,j}_k}{\rho^{i,j'}_k}+\frac{\rho^{i,j'}_k}{\rho^{i,j}_k}=+\infty\\[5mm]
\ds\mbox{ or }\lim_{k\rightarrow +\infty}\frac{|x_k^{i,j}-x_k^{i,j'}|}{\rho^{i,j}_k+\rho^{i,j'}_k}=+\infty\quad.
\end{array}
\rg.
\ee
Moreover,
\be
\label{III.11}
\forall\: i\in\{1,\ldots, n\}\quad\forall\: j\in\{1,\ldots, Q^i\}\quad\int_{B_{\rho^{i,j}_k}(x^{i,j}_k)}|d\vec{n}_{\vec{\xi}_k}|^2_{g_k}\ dvol_{g_k}>\ep_0\quad.
\ee
The set of balls $B_{\rho_k^{i,j}}$ are called ``bubbles'' associated to the sequence $\vec{\xi}_k\,$. 

\medskip

\noindent We have furthermore that for any $i\in\{1,\ldots, n\}$ and for any $j\in\{1,\ldots, Q^i\}$, the set of indices\footnote{this set might be empty.}
\[
I^{i,j}:=\lf\{j'\quad ; \quad x_k^{i,j'}\in B_{\rho^{i,j}_k}(x^{i,j}_k)\quad;\quad\frac{{\rho}^{i,j}_k}{\rho^{i,j'}_k}\rightarrow +\infty\rg\}
\]
is independent of $k$. It is called the set of ``bubbles contained in the bubble $B_{\rho^{i,j}_k}(x^{i,j})$''.

\medskip

\noindent For any $\al<1$ we denote 
\[
\forall\: i\in\{1,\ldots, n\}\quad\quad\Om_k^i(\al):=B_{\al}(a^i)\setminus \bigcup_{j=1}^{Q^i} B_{\al^{-1}\rho^{i,j}_k}(x^{i,j}_k)\quad,
\]
and 
\[
\forall\: i\in\{1,\ldots, n\}\quad\forall\: j\in\{1,\ldots, Q^i\}\quad\Om_k^{i,j}(\al):=\bigcup_{j'\in I^{i,j}}\lf(B_{\al\rho_k^{i,j}}(x^{i,j'}_k)\setminus \bigcup_{j''\in I^{i,j}} B_{\al^{-1}\rho_k^{i,j''}}(x_k^{i,j''})\rg)\quad.
\]
The sets $\Om_k^i(\al)$ and the $\Om_k^{i,j}(\al)$ are the ``$\al$-neck regions'' of $\vec{\xi}_k$. Let 
\[
\Om_k(\al):=\bigcup_{i=1}^n\lf(\Om_k^i(\al)\cup\,\bigcup_{j=1}^{Q^i}\Om_k^{i,j}(\al)\rg)
\]
Then there exists $0<\al<1$ independent of $k$, such that
\be
\label{III.12}
\begin{array}{l}
\ds\forall\: i\in\{1,\ldots, n\}\quad\forall\: j\in\{1,\ldots, Q^i\}\quad\mbox{ and}\:\:\:\forall\:\rho>0\quad\mbox{with }\:\:\:B_{2\rho}(x_k^{i,j})\setminus B_\rho(x)\subset \Om_k(\al)\\[5mm]
\ds\mbox{then }\quad\quad \int_{B_{2\rho}(x_k^{i,j})\setminus B_\rho(x)}|d\vec{n}_{\vec{\xi}_k}|^2_{g_k}\ dvol_{g_k}<\ep_0\quad,
\end{array}
\ee
and such that
\be
\label{III.13}
\begin{array}{l}
\ds\forall\: i\in\{1,\ldots, n\}\quad\forall\: j\in\{1,\ldots, Q^i\}\quad\quad\\[5mm]
\ds\mbox{then }\quad\quad\int_{B_{\rho_k^{i,j}}(x^{i,j}_k)\setminus \bigcup_{j'\in I^{i,j}}B_{\al\rho_k^{i,j}}(x^{i,j'}_k)}|d\vec{n}_{\vec{\xi}_k}|^2_{g_k}\ dvol_{g_k}>\ep_0\quad.
\end{array}
\ee
Given $0<\al<1$, we set
\[
\begin{array}{l}
\ds\forall\: i\in\{1,\ldots, n\}\quad\forall\: j\in\{1,\ldots, Q^i\}\quad\quad\\[5mm]
\ds r_k^{i,j}:=\inf\lf\{0<r<\rho^{i,j}_k\:\:\:;\ 
\begin{array}{l}
\ds B_\rho(x)\subset B_{\al^{-1}\rho_k^{i,j}}(x^{i,j}_k)\setminus \bigcup_{j'\in I^{i,j}} B_{\al\rho_k^{i,j}}(x^{i,j'}_k)\\[5mm]
\ds\int_{B_r(x)}\ |d\vec{n}_{\vec{\xi}_k}|^2_{g_k}\ dvol_{g_k}=\ep_0
\end{array}
\rg\}\quad.
\end{array}
\]
Then there holds
\be
\label{III.14}
\forall\: i\in\{1,\ldots, n\}\quad\forall\: j\in\{1,\ldots, Q^i\}\quad\quad\liminf_{k\rightarrow +\infty}\;\frac{r^{i,j}_k}{\rho^{i,j}_k}>0\quad.
\ee
\hfill $\Box$
\end{Prop}

We prove this proposition inductively by repeating finitely many iterations of the following lemma.
 
\begin{Lm}
\label{lm-III.1}{\bf [Energy tracking procedure]}
Let $\vec{\xi}_k$ be a sequence conformal weak immersions from the two-dimensional annulus $\Om_k:=B_{R_k}(0)\setminus B_{r_k}(0)\subset{\C}$ into ${\R}^m$,
with $R_k/r_k\rightarrow +\infty$. We assume that
\be
\label{IV.01}
\int_{B_{R_k}(0)\setminus B_{r_k}(0)}|\nabla \vec{n}_{\vec{\xi}_k}|^2\,dx\le\La<+\infty\quad.
\ee
Then for any $\ep_0>0$, one can extract a subsequence, still denoted $\vec{\xi}_k$, such that there exist a fixed  integer $N\in{\N}$ satisfying
\be
 \label{IV.1aaa}
 N\le 4\lf(\frac{\La}{\ep_0}+1\rg)\quad,
\ee
and a sequence of decreasing families of $(N+1)$ radii
 \[
 R_k^0=R_k>R_k^1>\ldots>R_k^N=r_k\quad,
 \]
 for which
 \be
\label{IV.2aaa}
\forall\: i\in \{0,\ldots, N-1\}\quad\quad\int_{\p B_{R_k^i}}|\nabla \vec{n}_{\vec{\xi}_k}|\ d{\mathcal H}^1< \sqrt{\frac{5\pi}{2}\ep_0}\quad .
\ee 
Moreover, the set $ I=\{0,\ldots, N-1\}$ can be decomposed into two disjoint subsets $I:=I_0\cup I_1$ with
\be
\label{IV.3aaa}
\forall\: i\in I_0\quad\quad\quad\log_2\frac{4}{3}<\tau_i:=\lim_{k\rightarrow +\infty}\log\frac{R_k^{i}}{R_k^{i+1}}< +\infty \quad,
\ee  
and
\be
\label{IV.2aaaaa}
\forall\: i\in I_1\quad\quad\quad\int_{B_{R_k^i}\setminus B_{R_k^{i+1}}}|\nabla \vec{n}_{\vec{\xi}_k}|^2\,dx<\ep\quad\quad\mbox{ and }\quad\quad\lim_{k\rightarrow +\infty}\log\frac{R_k^i}{R_k^{i+1}}=+\infty\quad.
\ee
\hfill $\Box$
\end{Lm}
{\bf Proof of lemma~\ref{lm-III.1}.}
We begin by constructing a decreasing sequence of radii $\{\rho_k^0=R_k>\rho_k^1>\ldots>\rho_k^{Q-1}>\rho_k^Q=r_k\}$ such that
\be
\label{IV.4}
\forall\: j\in\{0,\ldots, Q\}\quad\quad\int_{\p B_{\rho_k^j}(0)}|\nabla \vec{n}_{\vec{\xi}_k}|\ d{\mathcal H}^1<\sqrt{\frac{5\pi}{2}\ep_0}\quad
\ee
and
\be
\label{IV.4a}
\forall\: j\in\{0,\ldots, Q\}\quad\quad \log_2\frac{\rho_k^{j-1}}{\rho_k^j}>\log_2\frac{4}{3}\quad.
\ee
Moreover, for any $j\in J=\{0,\ldots, Q-1\}$, one of the following three possibilities holds
\begin{itemize}
\item[(i)]
\be
\label{IV.5}
\ep_0/2<\int_{B_{\rho_k^j}\setminus B_{\rho_k^{j+1}}}|\nabla \vec{n}_{\vec{\xi}_k}|^2\, dx<\ep_0
\ee
\item[(ii)]
\be
\label{IV.6}
\int_{B_{\rho_k^j}\setminus B_{\rho_k^{j+1}}}|\nabla \vec{n}_{\vec{\xi}_k}|^2\,dx<\ep_0/2
\ee
and either $\rho_k^{j+1}=\rho_k^{Q}=r_k$ or
\be
\label{IV.7}
\int_{B_{\rho_k^{j+1}}\setminus B_{\rho_k^{j+2}}}|\nabla \vec{n}_{\vec{\xi}_k}|^2\, dx>\ep_0/2\quad\quad\mbox{ and }\quad\quad\log_2\frac{\rho_k^{j+1}}{\rho_k^{j+2}}\le C
\ee
\item[(iii)]
\be
\label{IV.8}
\int_{B_{\rho_k^{j}}\setminus B_{\rho_k^{j+1}}}|\nabla \vec{n}_{\vec{\xi}_k}|^2\, dx>\ep_0/2\quad\quad\mbox{ and }\quad\quad\log_2\frac{\rho_k^{j}}{\rho_k^{j+1}}\le C\quad,
\ee
\end{itemize}
where\footnote{$[a]$ denotes the integer part of $a\in\R$, namely, the largest integer less than or equal to $a$.} $C:=\lf[2\La/\ep_0\rg]+2$. 

\medskip
\noindent
The construction of $\rho_k^j$ is done by induction on $j$. We take first $\rho_k^0:=R_k$. Assume that we have found some $\rho_k^j$ satisfying (\ref{IV.4})-(\ref{IV.4a}) and one of the alternatives (i)-(iii) for $l\le j$, and such that (\ref{IV.6}) fails for $j$. To construct $\rho_k^{j+1}$, we proceed as follows.
Let $\rho\ge r_k$ be such that 
\[
\int_{B_{\rho_k^j}\setminus B_\rho(0)}|\nabla \vec{n}_{\vec{\xi}_k}|^2\, dx=\ep_0/2\quad.
\]
If no such $\rho$ exists, we simply choose $\rho_k^{j+1}=r_k$.\\ 
For any $0<\al<1$  and $ R_k>t> \al^{-1}\,r_k$, there holds
\be
\label{IV.9}
\int_{\alpha t}^t ds\,\int_{\p B_s(0)}|\nabla\vec{n}_{\vec{\xi}_k}|\le t\  \sqrt{\pi(1-\al^2)}\ \sqrt{\int_{B_t\setminus B_{\al t}}|\nabla \vec{n}_{\vec{\xi}_k}|^2}\quad.
\ee
Suppose first that 
$$\int_{B_\rho\setminus B_{\rho/2}}|\nabla \vec{n}_{\vec{\xi}_k}|^2\le \ep_0/2\quad.$$ 
Applying (\ref{IV.9}) with $t=3\rho/4$ and $\al=2/3$, and using the mean value theorem, yields some $s\in [\rho/2,3\rho/4]$ such that
\[
\int_{\p B_{s}(0)}|\nabla \vec{n}_{\vec{\xi}_k}| d{\mathcal H}^1<\sqrt{\frac{5\pi}{2}\ep_0}\quad.
\]
Then $\rho_k^{j+1}:=s$ satisfies (\ref{IV.4}), (\ref{IV.4a}) and the alternative (\ref{IV.5}).
 
 \medskip
 
 \noindent 
Suppose next that
\[
\int_{B_\rho\setminus B_{\rho/2}}|\nabla \vec{n}_{\vec{\xi}_k}|^2> \ep_0/2\quad.
\]
Then there must exist some $l\in\{1,\ldots,\lf[2\La/\ep_0\rg]+1\}$ for which
\[
\int_{B_{2^{-l}\rho}\setminus B_{2^{-l-1}\rho}}|\nabla \vec{n}_{\vec{\xi}_k}|^2\le\ep_0/2\quad.
\]
Calling again upon (\ref{IV.9}) with $t=2^{-l}\rho$ and $\al=1/2$, and using the mean value theorem, we obtain some $s\in[2^{-l-1},2^{-l}]$ such that
\be
\label{IV.10}
\int_{\p B_{s}(0)}|\nabla \vec{n}_{\vec{\xi}_k}| d{\mathcal H}^1<\sqrt{\frac{3\pi}{2}\ep_0}\quad.
\ee
If $\rho_k^j\le 2\rho$, then we let $\rho_k^{j+1}:=s$, and note that (\ref{IV.4}), (\ref{IV.4a}), and the third alternative (\ref{IV.8}) with $C=\lf[2\La/\ep_0\rg]+2$ are satisfied.
If instead $\rho_k^j>2\rho$, we pick some $\sigma\in[2\rho,\rho]$ such that 
\[
\int_{\p B_{\sigma}(0)}|\nabla \vec{n}_{\vec{\xi}_k}| d{\mathcal H}^1<\sqrt{\frac{3\pi}{2}\ep_0}\quad,
\]
and set $\rho_k^{j+1}:=\sigma$ and $\rho_k^{j+2}:=s$. One verifies that (\ref{IV.4}), (\ref{IV.4a}), and the second alternative (\ref{IV.6})-(\ref{IV.7}) is satisfied.\\

\noindent
The above inductive construction must stop after a number of steps smaller than $[4\La/\ep_0]+2$. Indeed, at least every other index, an amount of $\ep_0/2$ of the (finite) Dirichlet
energy of $\vec{n}_{\vec{\xi}_k}$ is exhausted.

\medskip

To complete the proof of the lemma, we first extract a subsequence $\vec{\xi}_{k'}$ such that
\begin{itemize}
\item[-] $Q$ is independent of $k'$.
\item[-] for each $j\in\{0,\ldots, Q-1\}$, we have 
\[
\tau_j:=\lim_{k'\rightarrow +\infty}\log_2\frac{\rho_{k'}^{j}}{\rho_{k'}^{j+1}}\in{\R}\cup\{+\infty\}\quad.
\]
\end{itemize}  
We keep denoting this subsequence by $\vec{\xi}_k$. From the sequence $\rho_k^j$, we extract the desired decreasing sequence of radii $R_k^j$ by merging successive annuli such that
\[
\lim_{k'\rightarrow +\infty}\log_2\frac{\rho_{k'}^{j}}{\rho_{k'}^{j+1}}=\tau_j\in{\R}\quad
\]
and dropping the intermediate $\rho_k^j$. One easily verifies that the obtained sequence $R_k^j$ satisfies (\ref{IV.1aaa})-(\ref{IV.2aaaaa}), thereby completing the proof of lemma~\ref{lm-III.1}.
\hfill $\Box$

\medskip

\noindent{\bf Proof of proposition~\ref{pr-III.1}.}

Let $0<\ep_0<8\pi/3$. Consider
\[
\rho_k:=\inf\lf\{\rho\ ;\ \int_{B_\rho(x)}|d\vec{n}_{\vec{\xi}_k}|^2_{g_k}\ dvol_{g_k}=\ep_0\quad\forall\: x\in\Sigma\rg\}\quad,
\]
where $B_\rho(x)$ denotes the geodesic ball of center $x$ and radius $\rho$ for the constant scalar curvature $h_k$.
If 
\[
\limsup_{k\rightarrow +\infty}\rho_k>0\quad,
\]
we can extract a subsequence such that $\rho_k$ converges to a positive constant. Then there is no $a_i$. This would indeed require the concentration
of at least $8\pi/3$ of the energy, which is precluded from the result proved in \cite{Ri3}.  There is no bubble and the procedure stops.\\[1ex]
Alternatively, we consider the case when
\[
\lim_{k\rightarrow +\infty}\rho_k=0\quad.
\]
Let $x_k\in\Sigma$ be such that 
\[
\int_{B_{\rho_k}(x_k)}|d\vec{n}_{\vec{\xi}_k}|^2_{g_k}\ dvol_{g_k}=\ep_0\quad.
\]
Modulo the extraction of a subsequence, we have that $x_k$ converges to one of the points $a^i$. We set $x_k^{i,1}:=x_k$ and $\rho_k^{i,1}:=\rho_k$. By hypothesis, we may choose converging conformal coordinates around $a^i$. For notational convenience and clarity, since the converging metrics $h_k$ are uniformly equivalent to the flat metric, we shall work with flat geodesic balls. We can choose a subsequence $\vec{\xi}_k$ and a fixed radius $\al>0$ such that, for any 
\[
\lim_{k\rightarrow +\infty}\,\sup_{0<r\le\al}\lf\{r\:\; ;\quad\int_{B_\al(0)\setminus B_r(0)}|\nabla\vec{n}_{\vec{\xi}_k}|^2\,dx=\ep_0\rg\}=0
\]
and we apply lemma~\ref{lm-III.1} to this subsequence, still denoted $\vec{\xi}_k$, and to $R_k:=\al$ and $r_k:=\rho_k^{i,1}=\rho_k$. This gives then the existence of this family of radii
\[
R^0_k=\al>R^1_k>\ldots>R_k^N=\rho_k^{i,0}
\]
satisfying (\ref{IV.2aaa})-(\ref{IV.2aaaaa}) and 
\be
\label{III.16}
\lim_{k\rightarrow +\infty}\log\frac{R_k^0}{R_k^1}=+\infty\quad.
\ee
We merge successive radii $R^{i}_k$ and $R^{i+1}_k$ for $i\in I_1$ or for $i\in I_0$ such that
\[
\int_{B_{R^i_k}(x_k^{i,1})\setminus B_{R^{i+1}_k}(x_k^{i,1})}|\nabla\vec{n}_{\vec{\xi}}|^2\ dx_1\,dx_2\le\ep_0
\]
We then get a subsequence of radii $R^0_k=\al>R^{i_1}_k>,\ldots,>R_k^{i_P}=\rho_k^{i,0}$ for $P\le N$ such that either, for $i_l\in I_0$
\[
\lim_{k\rightarrow +\infty}\log\frac{R_k^{i_l}}{R_k^{i_{l+1}}}<+\infty\quad\mbox{ and }\quad\int_{B_{R^{i_l}_k}(0)\setminus B_{R^{i_{l+1}}_k}(0)}|\nabla\vec{n}_{\vec{\xi}}|^2\ dx_1\,dx_2\ge\ep_0\quad,
\]
or, when $i_l\in I_1$, we have
\[
\lim_{k\rightarrow +\infty}\log\frac{R_k^{i_l}}{R_k^{i_{l+1}}}=+\infty\quad\mbox{ and  }\quad\forall\:\rho\in(R_k^{i_{l+1}},R_k^{i_l}/2)\quad\int_{B_{2\rho}(x_k^{i,1})\setminus B_{\rho}(x_k^{i,1})}|\nabla\vec{n}_{\vec{\xi}}|^2\ dx_1\,dx_2<\ep_0
\]
We consider the smallest annulus $\Om^{i_l}_k:=B_{R_k^{i_l}}\setminus B_{R_k^{i_{l+1}}}$ of the first type $i_l\in I_0$. For such an $i_l$ we introduce
\[
\rho^{i_l}_k:=\inf\lf\{\rho<R^{i_l+1}_k\ ;\ \int_{B_\rho(x)}|d\vec{n}_{\vec{\xi}_k}|^2_{g_k}\ dvol_{g_k}=\ep_0\quad\mbox{ where }\quad x\in \Om^{i_l}_k \rg\}
\]
We now consider the following alternative. Either
\[
\limsup_{k\rightarrow +\infty}\frac{\rho^{i_l}_k}{R^{i_l}_k}>0\quad,
\]
we then extract a subsequence still denoted $\vec{\xi}_k$ such that 
\[
\lim_{k\rightarrow +\infty}\frac{\rho^{i_l}_k}{R^{i_l}_k}>0\quad,
\]
and we pass to the next $\Om_k^{i_{l'}}$, if there is any, where $i_{l'}\in I_0$, or we have
\[
\lim_{k\rightarrow +\infty}\frac{\rho^{i_l}_k}{R^{i_l}_k}=0\quad.
\]
In such a case we apply again the {\it energy tracking lemma~\ref{lm-III.1}} on the annulus $B_{R^{i_l}_k/2}(x)\setminus B_{\rho^{i_l}_k}(x)$ where $x\in \Om^{i_l}_k$ satisfies
\[
\int_{B_{\rho_k^{i_l}}(x)}|d\vec{n}_{\vec{\xi}_k}|^2_{g_k}\ dvol_{g_k}=\ep_0
\]
and we set $x_k^{i,2}:=x$ and $\rho_k^{i,2}:=\rho_k^{i_l}$. We keep this procedure running until all annuli of the type $I_0$ have been explored, and until the neighborhood of each $a^i$ has been considered. This can be done in a finite number of steps, since each application of lemma~\ref{lm-III.1} takes place around a ''bubble'' exhausting at least $\ep_0$ of the total energy of the sequence $\vec{\xi}_k$. With moderate effort, one checks that the sequence of finite families of bubbles constructed at the end of the procedure fulfills the desired requirements of proposition~\ref{pr-III.1}.\hfill $\Box$

\section{Construction of a Coulomb moving frame in each component of the neck region}
\reset

Let $\vec{\xi}_k$ be a family of weakly converging Willmore immersions of $\Sigma$ fulfilling the conclusions of lemma~\ref{lm-I.1}.
 For some $\ep_0$ chosen smaller than the constant $\ep$ appearing in the $\ep$-regularity theorem I.5 in \cite{Ri2}, we apply proposition~\ref{pr-III.1} so as to obtain a subsequence for
which we can decompose the surface into ``converging regions" (complements of small neighborhoods of the points $a^i$), bubbles, and neck regions. In each annulus $B_{R_k}(x_k^{i,j})\setminus B_{r_k}(x_k^{i,j})\subset \Om_k(\al)$ of a part of the neck region $\Om_k(\al)$ satellite to some bubble centered on $x_k^{i,j}$, the condition (\ref{III.12}) states that
 \[
 \quad\quad \sup_{r_k<\rho<R_k}\int_{B_{2\rho}(x_k^{i,j})\setminus B_\rho(x_k^{i,j})}|\nabla\vec{n}_{\vec{\xi}_k}|\,dx<\ep_0\quad.
 \]
 Combining this fact with the $\ep$-regularity theorem I.5 in \cite{Ri2} gives that 
 \be
 \label{V.1}
 \forall\: x\in B_{R_k/2}(x_k^{i,j})\setminus B_{2r_k}(x_k^{i,j})\quad\quad|\nabla\vec{n}_{\vec{\xi}_k}|\le C_m\ \frac{\sqrt{\ep_0}}{|x-x_k^{i,j}|}\quad.
 \ee
 This pointwise estimate translates into a control of the $L^2$-weak norm of $\nabla\vec{n}_{\vec{\xi}}$ in the annulus, namely:
 \be
 \label{V.2}
 \|\nabla\vec{n}_{\vec{\xi}_k}\|_{L^{2,\infty}(B_{R_k/2}(x_k^{i,j})\setminus B_{2r_k}(x_k^{i,j}))}\le C_m\,\sqrt{\ep_0}\quad.
 \ee
This $L^2$-weak estimate for the gradient of the Gauss map in neck regions already appears in previous works \cite{LiRi1}, \cite{LiRi2} and \cite{Ri5}, where it plays a decisive role. It is also of utmost importance in the present paper, as it will enable us in lemma~\ref{lm-AA.1} to construct a $W^{1,2}$ Coulomb gauge without having to resort to the stronger $L^2$-control of the gradient of the Gauss map commonly required in similar situations. Our results are generic in nature, and they should be helpful to the study of other gauge invariant problems.

\smallskip

We first need an extension lemma.

\begin{Lma}
\label{lm-A.2}
Let $\vec{n}\in W^{1,2}(B_2(0)\setminus B_1(0), G_{m-2}({\R}^m))$ where $G_{m-2}({\R}^m)$ denotes the Grassmannian manifold of oriented
$(m-2)$-planes in ${\R}^m$. There exists $\delta(m)>0$ depending only  on $m$ such that if
\be
\label{A.14}
\int_{\p B_1(0)}|\nabla\vec{n}|\ d\theta\le \delta(m)\quad,
\ee
then we can extend $\vec{n}$ to $\hat{\vec{n}}\in W^{1,2}(B_2(0), G_{m-2}({\R}^m))$ with $\hat{\vec{n}}=\vec{n}$ on the annulus $B_2(0)\setminus B_1(0)$, 
and furthermore
\be
\label{A.15}
\int_{B_2(0)}|\nabla\hat{\vec{n}}|^2\,dx\le  C_m\ \int_{B_2(0)\setminus B_1(0)}|\nabla\vec{n}|^2\, dx\quad,
\ee
where $C_m>0$ depends only on $m$. 
\hfill $\Box$
\end{Lma}
{\bf Proof of lemma~\ref{lm-A.2}.}
Suppose that
\be
\label{A.14bis}
\int_{\p B_1(0)}|\nabla\vec{n}|\ d\theta\le \delta\quad,
\ee
for some $\delta>0$. Then
\[
\|\vec{n}-\vec{n}(1,0)\|_{L^\infty(\p B_1(0))}\le\int_{\p B_1(0)}\lf|\frac{\p\vec{n}}{\p\theta}\rg|\ d\theta\le \delta\quad.
\]
Thus, if $\delta>0$ chosen small enough and independently of $\vec{n}(1,0)$, the plane $\vec{n}(p)$ is contained in a small geodesic ball $B_\delta(\vec{n}(1,0))\subset G_{m-2}({\R}^m)$ on which coordinates 
$y=(y_1,\ldots, y_N)$ exist. We may further choose these coordinates in such a way that $B_\delta(\vec{n}(1,0))$ corresponds to the ball $B^N_\delta(0)$ of radius $\delta$ and centered on $0\in{\R}^N$. Owing to the compactness of $G_{m-2}({\R}^m)$, we may always arrange for $\|\nabla y^{-1}\|_{L^\infty(B_\delta(0))}$ to depend only on $m$.   We have\footnote{On the left-hand side of (\ref{A.16}), we use the isometric embedding of $G_{m-2}({\R}^m)$ into some 
${\R}^{K_m}$, and $\vec{n}$ is viewed as an $\R^{K_m}$-valued map.}
\be
\label{A.16}
\|\vec{n}\|_{H^{1/2}(\p B_1(0))}\simeq\sum_{l=1}^N\|{n}_l\|_{H^{1/2}(\p B_1(0))}\quad,
\ee
where ${n}_l$ are the coordinates of $\vec{n}$ in the system $\{y_l\}_{l=1,\ldots,N}$. Let $\hat{n}_l$ be the harmonic extension of $n_l$ into $B_1(0)$. Since the $\R^N$-valued
map $(n_l)$ is included in the ball $B^N_\delta(0)$, the maximum principle yields that $(\hat{n}_l)_{l=1,\ldots, N}\in B^N_\delta(0)$. Moreover,
\be
\label{A.17}
\int_{B_1(0)}|\nabla\hat{n}_l|^2\ \le C\sum_{l=1}^N\|{n}_l\|^2_{H^{1/2}(\p B_1(0))}\le c_m\ \|\vec{n}\|^2_{H^{1/2}(\p B_1(0))}\le  C_m\ \int_{B_2(0)\setminus B_1(0)}|\nabla{\vec{n}}|^2\quad,
\ee
where the constants $c_m$ and $C_m$ depend only on $m$. Note that in the last inequality, we have again used the characterization of $\vec{n}$  as an $\R^{K_m}$-valued map. We set $\hat{\vec{n}}:=y^{-1}((\hat{n}_l))$. Since $\|\nabla y^{-1}\|_{L^\infty(B_\delta(0))}$
only depends on $m$, (\ref{A.17}) implies the desired (\ref{A.15}), thereby proving lemma~\ref{lm-A.2} with $\delta(m)=\delta$ small enough.\hfill$\Box$

\medskip

The main result of this section also requires an integrability-by-compensation lemma for second-order equations in divergence form. It was originally proved by Y.~Ge in \cite{Ge} through a different method from the one used here. 

\begin{Lma}\label{weakwente}{\bf [Wente estimates in Lorentz $L^{p,q}$ spaces]}
Consider the divergence-form problem
\be\label{l0l}
\left\{\begin{array}{rcll}\Delta\varphi&=&\nabla^\perp a\cdot\nabla b&,\:\:\:{in}\:\:\:B_1(0)\\[1ex]
\varphi&=&0&,\:\:\:{on}\:\:\:\partial B_1(0)\quad,
\end{array}\right.
\ee
where $\nabla a\in L^{2,\infty}(B_1(0))$, and $\nabla b\in L^{p,q}(B_1(0))$ for some $p\in(1,\infty)$ and $q\in[1,\infty]$. Then there holds
\be\label{ww0}
\Vert\nabla\varphi\Vert_{L^{p,q}(B_1(0))}\;\leq\;C_{p,q}\,\Vert\nabla a\Vert_{L^{2,\infty}(B_1(0))}\Vert\nabla b\Vert_{L^{p,q}(B_1(0))}\quad,
\ee
for some constant $C_{p,q}>0$.\hfill $\Box$
\end{Lma}
{\bf Proof of lemma~\ref{weakwente}.} Owing to the interpolative nature of the Lorentz space $L^{p,q}$, it suffices to show that (\ref{ww0}) holds for $p=q\in(1,\infty)$. Consider first the case when $p>2$. According to the H\"older inequality for Lorentz spaces, there holds
\be
L^{2,\infty}\cdot L^{p}\;\subset\;L^{\frac{2p}{2+p}\,,p}\quad.
\ee
Since $p>2$, we call upon the usual Calderon-Zygmund theorem and the Sobolev embedding theorem to obtain that $\nabla\varphi\in W^{1,\frac{2p}{2+p}}\subset L^{p}\,$ with the estimate (\ref{ww0}), as claimed. \\
Next, when $1<p\le 2$, we use the divergence-form structure of the equation. Note that
\be\label{ww1}
\Delta\varphi\;=\;div\,(b\,\nabla^\perp a)\quad,
\ee
and that
\be\label{ww2}
b\,\nabla^\perp a\;\in\;W^{1,p}\cdot L^{2,\infty}\;\subset\; L^{\frac{2p}{2-p},p}\cdot L^{2,\infty}\;\subset\;L^p\quad,
\ee
where we have again used the Sobolev embedding theorem and the H\"older inequality for Lorentz spaces. The desired (\ref{ww0}) follows immediately from (\ref{ww1}) and (\ref{ww2})
and this concludes the proof of lemma~\ref{weakwente}.\hfill $\Box$  

\medskip

We next state and prove the main result of this section. 
\medskip

\begin{Lma}
\label{lm-AA.1}{\bf[$W^{1,2}$-controlled Coulomb frame from small $L^{2,\infty}$-control of the second fundamental form]}
Let ${\bn}\in W^{1,2}(B_1(0)\setminus B_r(0), G_{m-2}({\R}^m))$ where $0\le r<1/2$. There exists $\eta(m)>0$ depending only on $m$ (but not on $r$) such that if
\be
\label{bA.14}
\Vert\nabla{\bn}\Vert_{L^{2,\infty}(B_1(0))}\;\le\;\eta(m)\:,
\ee
and if
\be
\label{bA.15}
\int_{\p B_r(0)}|\nabla \vec{n}|\ dl_{\p B_r}\le\;\eta(m)\quad,
\ee
where $dl_{\p B_r}$ is the length form $r\,d\theta$ on $\p B_r$,
then one can construct an associated Coulomb frame $\{{\bbe_1},{\bbe_2}\}$ satisfying
\[
\star\,{\bn}\,=\,{\bbe_1}\wedge{\bbe_2}\qquad\mbox{and}\qquad div({\bbe_1}\cdot\nabla{\bbe_2})\,=\,0\quad.
\]
Furthermore,
\be
\label{bbA.15}
\int_{B_1(0)\setminus B_r(0)}|\nabla\bbe_1|^2+|\nabla\bbe_2|^2\,\leq\,C\,\int_{B_1(0)\setminus B_r(0)}|\nabla{\bn}|^2
\ee
and
\be
\label{bbA.16}
\sum_{i=1}^2\|\nabla\vec{e}_i\|_{L^{2,\infty}(B_1(0)\setminus B_r(0))}\le C\,\|\nabla\vec{n}\|_{L^{2,\infty}(B_1(0)\setminus B_r(0))}\quad,
\ee
for some constant $C$ depending on $m$, but neither on $r$ nor on $\vec{n}$.\hfill $\Box$
\end{Lma}
{\bf Proof of lemma~\ref{lm-AA.1}.}
Owing to lemma~\ref{lm-A.2}, we can always reduce to the case when $r=0$. Indeed, we have seen that the extension $\hat{\vec{n}}$ in $B_r(0)$ of the restriction
of $\vec{n}$ to $B_{2r}(0)\setminus B_{r}(0)$ satisfies
\[
\int_{B_{r}(0)}|\nabla\hat{\vec{n}}|^2\, dx\le C_m\,\int_{B_{2r}(0)\setminus B_{r}(0)}|\nabla{\vec{n}}|^2\, dx\le C_m\ \eta(m)\quad,
\]
In particular, if lemma~\ref{lm-A.2} holds for $r=0$ with the constant $\eta(m)>0$, then it also holds for an arbitrary $0<r<1/2$ with the constant $\eta(m)/(C_m+1)$ playing the role of $\eta(m)$ in the statement. We shall thus focus on the case $r=0$. 

\medskip

We first suppose that $\vec{n}$ is smooth.
According to Theorem 5.2.1 in \cite{Hel}, there exists a smooth moving orthonormal two-frame $\{{\bbf_1},{\bbf_2}\}$ in $W^{1,2}(B_1(0),\R^m)$ with $\,\star\,{\bn}={\bbf_1}\wedge{\bbf_2}$. Its energy is however not controlled, since we do not assume that $\Vert\nabla{\bn}\Vert_{L^{2}(B_1(0))}$ is smaller than the required threshold $8\pi/3$. \\
For each $r\in[0,1]$, let $\,\{{\bbf}_{1,r}(x),{\bbf}_{2,r}(x)\}:=\{{\bbf_{1}}(rx),{\bbf_{2}}(rx)\}$. We minimize
\be
F_r(\theta_r)\;:=\;\int_{B_r(0)}\big|\nabla\theta_r+({\bbf}_{1,r}\cdot\nabla{\bbf}_{2,r})\big|^2
\ee
over all rotations $\theta_r\in W^{1,2}(B_1(0),\R)$. As explained in Theorem 4.1.1 from \cite{Hel}, for each $r$, the minimum of $F_r$ is attained at some frame $\{\bbe_{1,r}\,,\bbe_{2,r}\}$ satisfying
\be\label{lar-1}
\bbe_{1,r}+i\,\bbe_{2,r}\,=\,{e}^{i\theta_r}\big({\bbf}_{1,r}+i\,{\bbf}_{2,r}\big)
\ee
 and
 \be
 \label{larn-1}
 \lf\{
 \begin{array}{l}
div\,(\bbe_{1,r}\cdot\nabla\bbe_{2,r})\,=\,0\quad\mbox{ on}\:\:\:B_r(0)\\[5mm]
\vec{e}_{1,r}\cdot\p_\nu\vec{e}_{2,r}\,=\, 0\quad\mbox{ on}\:\:\:\p B_r(0)\quad.
\end{array}
\rg.
\ee
Thus, there exists $\la_r$  equal to zero on $\p B_r(0)$ and solving
\[
\nabla\la_r=-\vec{e}_{1,r}\cdot\nabla^\perp\vec{e}_{2,r}\quad.
\]
It satisfies in particular
\be\label{lar0}
\left\{\begin{array}{rcll}\Delta\lambda_r&=&\nabla^\perp\bbe_{1,r}\cdot\nabla\bbe_{2,r}&,\:\:\:{in}\:\:\:B2_r(0)\\[1ex]
\lambda_r&=&0&,\:\:\:{on}\:\:\:\partial B_r(0)\quad.
\end{array}\right.
\ee
Using the Wente inequality lemma~\ref{weakwente} with $(p,q)=(2,\infty)$ yields
\be\label{lar1}
\Vert\nabla\la_r\Vert_{L^{2,\infty}(B_r(0))}\,\le\;C_0\,\Vert\nabla\bbe_{1,r}\Vert_{L^{2,\infty}(B_r(0))}\Vert\nabla\bbe_{2,r}\Vert_{L^{2,\infty}(B_r(0))}\quad,
\ee
for some constant $C_0$ independent of $r$. Moreover, as shown in the proof of Lemma 5.1.4 in \cite{Hel}, an elementary computation gives
\be\label{lar2}
|\nabla\bbe_{1,r}|^2+\,|\nabla\bbe_{2,r}|^2\;=\;2\,|\nabla\la_r|^2+\,|\nabla{\bn}|^2\quad.
\ee
Combining altogether (\ref{lar1}) and (\ref{lar2}) gives the inequality
\be\label{lar3}
\Vert\nabla\la_r\Vert^2_{L^{2,\infty}(B_r(0))}\,+\,\frac{1}{4C_0}\big(C_0\Vert\nabla{\bn}\Vert_{L^{2,\infty}}-1\big)\Vert\nabla\la_r\Vert_{L^{2,\infty}(B_r(0))}\,+\,C_0\Vert\nabla{\bn}\Vert_{L^{2,\infty}(B_r(0))}\;\geq\;0\quad.
\ee
Accordingly, there exists some small enough threshold $\eta(C_0)\equiv\eta(m)$ for which the hypothesis that
$\,\Vert\nabla{\bn}\Vert_{L^{2,\infty}(B_1(0))}<\eta(m)\,$ guarantees that the range of $\,\Vert\nabla\la_r\Vert_{L^{2,\infty}(B_r(0))}$ is of the form $[0,\al(\eta)]\cup[\beta(\eta),\infty)$ for some $0<\alpha(\eta)<\beta(\eta)$. 
One easily checks from (\ref{lar3}) that $\alpha(\eta)\leq C_1\|\nabla{\bn}\Vert_{L^{2,\infty}(B_1(0))}$ with the constant $C_1$ depending only on $C_0$, and thus only on $m$.

\medskip
\noindent
As ${\bbf}_{j,r}(x)={\bbf}_{j}(rx)$ and ${\bbf}_{j}$ belongs to $W^{1,2}(B_1(0))$, it is clear that ${\bbf}_{j,r}$ is continuous in the parameter $r$. From this and the definition of $\theta_r$, it follows that $\theta_r$ is likewise continuous in $r$, and hence from (\ref{lar-1}) and (\ref{lar0}), that $\bbe_{j,r}$ and $\la_r$ are continuous in $r$. We see that $\Vert\nabla\la_r\Vert_{L^{2,\infty}(B_r(0))}$ takes the value 0 at $r=0$. From the above discussion, we deduce that at $r=1$, 
\be\label{lar4}
\Vert\nabla\la\Vert_{L^{2,\infty}(B_1(0))}\;\equiv\;\Vert\nabla\la_1\Vert_{L^{2,\infty}(B_1(0))}\;\leq\;\alpha(\eta)\;\le\;C_1\,\|\nabla{\bn}\Vert_{L^{2,\infty}(B_1(0))}\quad.
\ee
Let now $\{\bbe_1,\bbe_2\}:=\{\bbe_{1,1},\bbe_{2,1}\}$. By definition, this Coulomb frame satisfies 
\be
\star\,{\bn}\,=\,{\bbe_1}\wedge{\bbe_2}\quad,\qquad div({\bbe_1}\cdot\nabla{\bbe_2})\,=\,0\quad,\qquad\nabla^\perp\la\,=\,\bbe_1\cdot\nabla\bbe_2\,=\,-\,\bbe_2\cdot\nabla\bbe_1\quad.
\ee
Note next that
\be\label{lar5}
|\nabla\bbe_j|^2\;=\;|\bbe_{j'}\cdot\nabla\bbe_j|^2\,+\,|\pi_{{\bn}}\nabla\bbe_j|^2\;\le\;|\nabla\la|^2\,+\,|\nabla{\bn}|^2\:,\qquad \mbox{for}\quad (j,j')\in\big\{(1,2)\,,(2,1)\big\}\quad.
\ee
In particular using (\ref{lar4}), we find
\be\label{lar6}
\Vert\nabla\bbe_j\Vert_{L^{2,\infty}(B_1(0))}\;\leq\;\Vert\nabla\la\Vert_{L^{2,\infty}(B_1(0))}+\Vert\nabla{\bn}\Vert_{L^{2,\infty}(B_1(0))}\;\le\;(C_1+1)\,\|\nabla{\bn}\Vert_{L^{2,\infty}(B_1(0))}\quad,
\ee
thereby giving the required (\ref{bbA.16}). \\
Calling again upon the Wente inequality of lemma~\ref{weakwente} with $(p,q)=(2,2)$, we obtain from (\ref{lar0}) with $r=1$ that
\begin{eqnarray*}
\Vert\nabla\la\Vert_{L(B_1(0))}&\le&C_2\,\Vert\nabla\bbe_1\Vert_{L^{2,\infty}(B_1(0))}\Vert\nabla\bbe_2\Vert_{L(B_1(0))}\\[1ex]
&\le&C_2\,(C_1+1)\,\eta(m)\,\big(\Vert\nabla\la\Vert_{L(B_1(0))}\,+\,\Vert\nabla{\bn}\Vert_{L(B_1(0))}\big)\quad,
\end{eqnarray*}
where we have used (\ref{lar5}), (\ref{lar6}), and the hypothesis (\ref{bA.14}). Here $C_2$ is a constant that depends only on $m$. The latter shows that $\,\Vert\nabla\la\Vert_{L(B_1(0))}\leq C_3\,\eta(m)\,\Vert\nabla{\bn}\Vert_{L(B_1(0))}$, which, once introduced into (\ref{lar5}) gives the announced
\be
\int_{B_1(0)}|\nabla\bbe_1|^2+|\nabla\bbe_2|^2\,\leq\,C\int_{B_1(0)}|\nabla{\bn}|^2\quad,
\ee
where the constant $C$ depends only on $m$.\\[1ex]
In the general case when ${\bn}$ is not smooth, we use the limiting process outlined in step 6 of the proof of Lemma 5.1.4 from \cite{Hel}, and based on the known fact \cite{ScU} that $C^\infty(B_1(0),G_{m-2}(\R^m))$ is dense into $W^{1,2}(B_1(0), G_{m-2}({\R}^m))$. This concludes the proof of lemma~\ref{lm-AA.1}.\hfill $\Box$

\medskip

\section{General estimates in the neck regions of weakly converging weak conformal immersions} 
\reset

\begin{Lma}
\label{lm-A.a00}{\bf [Estimate of the ``Euler characteristic production'' in neck regions]}
There exists a constant $\eta(m)>0$ with the following property. Let $\eta<\eta(m)$, let $R>4r>0$, and let $\vec{\xi}$ be a conformal weak immersion from $B_{R}(0)\setminus B_r(0)$ into ${\R}^m$ with $L^2$-bounded second fundamental form satisfying
\be
\label{A.b0}
\|\nabla\vec{n}_{\vec{\xi}}\|_{L^{2,\infty}(B_R\setminus B_r)}\le\eta\quad,
\ee
\be
\label{A.b00}
\int_{\p B_r(0)} |\nabla\vec{n}_{\vec{\xi}}|\ dl_{\p B_r(0)}\le\eta\quad,
\ee
and
\be
\label{A.b000}
\int_{B_R\setminus B_{R/2}\cup B_{2r}\setminus B_r}|\nabla\vec{n}_{\vec{\xi}}|^2\, dx\le\eta\quad.
\ee
Then there exists a  positive constant $C_m$ depending only on $m$, but not on the conformal type of the annulus $\log\frac{R}{r}$, such that
\be
\label{A.b02}
\lf|\int_{B_{R}(0)\setminus B_{r}(0)}K_{\vec{\xi}}\ dvol_{g_{\vec{\xi}}}\rg|\le C_{m}\,\eta\quad,
\ee
where $K_{\vec{\xi}}$ denotes the Gaussian curvature associated with the immersion $\vec{\xi}\,$. 
\hfill$\Box$
\end{Lma}
{\bf Proof of lemma~\ref{lm-A.a00}.}
Let $\Om=B_R(0)\setminus B_r(0)$. According to lemma~\ref{lm-AA.1} (applicable here owing to the hypotheses (\ref{A.b0}) and (\ref{A.b00})), there exists a Coulomb moving frame $$\{\vec{e}_{1},\vec{e}_{2}\}\in (W^{1,2}(\Om,S^{m-1}))^2$$ satisfying
$$
\star\,\vec{n}_{\vec{\xi}}=\vec{e}_{1}\wedge\vec{e}_{2}
$$
and
\be
\label{A.0a2}
\sum_{i=1}^2\|\nabla \vec{e}_{i}\|_{L^{2,\infty}(\Om)}\le C_1\, \|\nabla\vec{n}\|_{L^{2,\infty}(\Om)}\le C_2\, {\eta}\quad,
\ee
where the constants $C_j$ depend only on $m$. The mean value theorem yields the existence of some $\al\in(1/2,1)$ with
\begin{eqnarray}\label{A.ba3}
\ds\,\sum_{i=1}^2\,\int_{\p (B_{\al R}\setminus B_{\al^{-1} r})}|\nabla \vec{e}_{i}|&\le& C_3\, r\ \sum_{i=1}^2\|\nabla \vec{e}_{i}\|_{L^1(B_{2r}\setminus B_{r})}
+ C_4\,R^{}\ \sum_{i=1}^2\|\nabla \vec{e}_{i}\|_{L^1(B_{R}\setminus B_{R/2})}\nonumber\\[.5ex]
&\le&\ds C_5\ \sum_{i=1}^2\|\nabla \vec{e}_{i}\|^{}_{L^{2,\infty}(\Om)}\;\;\le\;\; C_5\,C_2\,{\eta}\quad,
\end{eqnarray}
where $C_5$ is a universal constant.\\
The Liouville equation states that
\be
\Delta \la\,=\,-\,e^{2\la}K_{\vec{\xi}}\quad,
\ee
where $\la$ is the conformal parameter of the conformal immersion $\vec{\xi}$. From the way the Coulomb frame $\{\bbe_1,\bbe_2\}$ was constructed, there also holds on the annulus $B_R(0)\setminus B_r(0)$\,:
\be
\Delta \la\,=\,-\,\nabla\bbe_1\cdot\nabla^\perp\bbe_2\,=\,-\,div\,\big(\bbe_1\cdot\nabla^\perp\bbe_2\big)\:.
\ee
Hence, we find using (\ref{A.ba3})
\be\label{A.ba35}
\begin{array}{l}
\ds\lf|\int_{B_{\al\, R}(0)\setminus B_{\al^{-1}\, r}(0)}K_{\vec{\xi}}\ e^{2\la}\, dx\rg|=\lf|\int_{B_{\al\, R}(0)\setminus B_{\al^{-1}\, r}(0)}\nabla\vec{e}_{1}\cdot\nabla^\perp\vec{e}_{2}\, dx\rg|\\[5mm]
\ds\quad\quad\le\int_{\p B_{\al\, R}(0)}\ |e_{1}\cdot\p_\tau e_{2}|\ dl+\int_{\p B_{\al^{-1}\, r}(0)}\ |e_{1}\cdot\p_\tau e_{2}|\ dl\le C\ \eta\quad,
\end{array}
\ee
where $\tau$ denotes the unit tangent vector to the circles $\p B_\rho(0)$, namely  $\tau=\rho^{-1}\ \p_\theta$.\\
Moreover, owing to the hypothesis (\ref{A.b000}) and the fact that $\al\in(1/2,1)$, we have that
\be
\label{A.ba4}
\int_{B_{\al^{-1}r}(0)\setminus B_{r}(0)}|K_{\vec{\xi}}|\ e^{2\la}\, dx\le 2^{-1}\int_{B_{\al^{-1}r}(0)\setminus B_{r}(0)}|\nabla\vec{n}_{\vec{\xi}}|^2\, dx\le
2^{-1}\,\eta\quad,
\ee
and similarly
\be
\label{A.ba5}
\int_{B_{R}(0)\setminus B_{\al R}(0)}|K_{\vec{\xi}}|\ e^{2\la}\, dx\le 2^{-1} \int_{B_{R}(0)\setminus B_{\al R}(0)}|\nabla\vec{n}_{\vec{\xi}}|^2\, dx\le
2^{-1}\,\eta\quad.
\ee
Altogether (\ref{A.ba35}), (\ref{A.ba4}), and (\ref{A.ba5}) yield the announced (\ref{A.b02}). \hfill $\Box$

\medskip

\begin{Lma}
\label{lm-A.0}{\bf[Uniform Harnack estimates for the conformal factor in neck regions]}
There exists a constant $\ep(m)>0$ with the following property. Let $1<R<\infty$. If $\vec{\xi}$ is a Lipschitz conformal immersion from $\Om_R:=B_{R}(0)\setminus B_{R^{-1}}(0)$ into ${\R}^m$ satisfying
\be
\label{A.a1}
\sup_{r\in (R^{-1},R/2)}\int_{B_{2\, r}(0)\setminus B_r(0)}|\nabla\vec{n}_{\vec{\xi}}|^2\, dx\le\ep(m)\quad,
\ee
and\footnote{as before, $\la$ is the conformal parameter of $\vec{\xi}$. Namely,  $|\p_{x_1}\vec{\xi}|={e}^{\la}=|\p_{x_2}\vec{\xi}|$.}
\be
\label{A.a2}
\|\nabla \la\|_{L^{2,\infty}(\Om_R)}=A<+\infty\quad.
\ee
Then for any $4R^{-2}<\al<1$ there exists $C_\al(A)>0$, depending only on $\al$ and $A$ but not on $R$, such that
\be
\label{A.a3}
\sup_{r\in (4\,R^{-1},\al R)} \frac{\sup_{x\in B_{\al^{-1} r}(0)\setminus B_r(0)}e^{\la(x)}}{\inf_{x\in B_{\al^{-1} r}(0)\setminus B_r(0)}e^{\la(x)}}\le C_\al\quad.
\ee
\hfill$\Box$
\end{Lma}
{\bf Proof of lemma~\ref{lm-A.0}.}
Let $r\in (4R^{-1},R/4)$. From (\ref{A.a1}), we have
\be
\label{A.a4}
\int_{B_{3\, r}(0)\setminus B_{r/4}(0)}|\nabla\vec{n}_{\vec{\xi}}|^2\, dx\le 4\ep(m)\quad.
\ee
With the help of the Fubini and mean value theorems, we deduce from (\ref{A.a1}) the existence of $\rho\in (r/4,r/2)$ such that
\be
\label{A.a5}
\int_{\p B_\rho(0)} |\nabla\vec{n}_{\vec{\xi}}|^2\,\le \frac{4}{r}\ep(m)\le\frac{2}{\rho}\ep(m)\quad,
\ee
hence
\be
\label{A.a6}
\int_{\p B_\rho(0)} |\nabla\vec{n}_{\vec{\xi}}|\, dx\le\sqrt{4\pi\, \ep(m)}\quad.
\ee
Following lemma \ref{lm-A.2}, we choose $\sqrt{4\pi\,\ep(m)}<\delta(m)$ so as to obtain an extension $\hat{\vec{n}}$ of $\vec{n}$ in $B_\rho(0)$ with
\be
\label{A.a7}
\int_{B_{4\, r}(0)}|\nabla\hat{\vec{n}}|^2\, dx\le  4\, C_m\,\ep(m)\quad,
\ee
where $C_m$ is the constant depending only on $m$ in the statement of lemma~\ref{lm-A.2}. Insuring further that $12\, C_m\, \ep_m\le 8\pi$,
we can apply lemma 5.1.4 of \cite{Hel} in order to find a framing $\{\vec{e}_1,\vec{e}_2\}\in (W^{1,2}(B_{4\, r}(0)))^2$ satisfying
 $\star\,\hat{\vec{n}}=\vec{e}_1\wedge\vec{e}_2$ and
\be
\label{A.a8}
\int_{B_{4\, r}(0)} \sum_{i=1}^2|\nabla\vec{e}_i|^2\, dx\le C\ \int_{B_{4\, r}(0)}|\nabla\hat{\vec{n}}|^2\, dx\le  4\,C\, C_m\, \ep(m)\quad.
\ee
The conformal parameter satisfies
\be
\label{A.a9}
\Delta\la=\nabla^\perp\vec{e}_1\cdot\nabla\vec{e}_2\quad\quad\mbox{ in }B_{4r}(0)\setminus B_{r/2}(0)\quad.
\ee
Let $\mu$ be the solution of
\be
\label{A.a10}
\lf\{
\begin{array}{l}
\ds\Delta\mu=\nabla^\perp\vec{e}_1\cdot\nabla\vec{e}_2\quad\quad\mbox{ in }B_{4r}(0)\setminus B_{r/2}(0)\\[5mm]
\ds\mu=0\quad\quad\quad\mbox{ on }\p(B_{4r}(0)\setminus B_{r/2}(0))\quad.
\end{array}
\rg.
\ee
A standard Wente estimate (see theorem 3.3.8 in \cite{Hel}) yields on one hand
\be
\label{A.a11}
\|\mu\|_{L^\infty(B_{4r}(0)\setminus B_{r/2}(0))}\le (2\pi)^{-1}\, \|\nabla\vec{e}_1\|_{L^2}\, \|\nabla\vec{e}_2\|_{L^2}\le C\, \ep(m)\quad,
\ee
and on the other hand
\be
\label{A.a12}
\|\nabla\mu\|_{L(B_{4r}(0)\setminus B_{r/2}(0))}\le \sqrt{3/32\pi}\, \|\nabla\vec{e}_1\|_{L^2}\, \|\nabla\vec{e}_2\|_{L^2}\le C\, \ep(m)\quad.
\ee
Let $\nu:=\la-\mu$ be the harmonic function in $B_{4r}(0)\setminus B_{r/2}(0)$ equal to $\la$ on the boundary. Owing to (\ref{A.a2}) and (\ref{A.a12}), there holds
\be
\label{A.a13}
\|\nabla\nu\|_{L^{2,\infty}(B_{4r}(0)\setminus B_{r/2}(0))}\le C\ \ep(m)+A\quad.
\ee
Let $\ov{\nu}$ be the average of $\nu$ on $B_{4r}(0)\setminus B_{r/2}(0))$. The Sobolev-Poincar\'e inequality gives
\begin{eqnarray}
\label{A.a14}
\ds r^{-2}\int_{B_{4r}(0)\setminus B_{r/2}(0)}|\nu-\ov{\nu}|\, dx&\le& C_1\, r^{-1}\int_{B_{4r}(0)\setminus B_{r/2}(0)}|\nabla\nu|\, dx\nonumber\\[5mm]
&\le&\ds C_2\,\|\nabla\nu\|_{L^{2,\infty}(B_{4r}(0)\setminus B_{r/2}(0))}\le C_3\ \ep(m)+C_2\,A\quad,
\end{eqnarray}
where $C_j$ are universal constants. Since $\nu-\ov{\nu}$ is harmonic, we deduce that on a smaller annulus the $L^\infty$-norm
of $(\nu-\ov{\nu})$ is controled by $\ep(m)$ and $A$\,:
\be
\label{A.a15}
\|\nu-\ov{\nu}\|_{L^\infty(B_{2r}(0)\setminus B_{r}(0)}\le C\  r^{-2}\int_{B_{4r}(0)\setminus B_{r/2}(0)}|\nu-\ov{\nu}|\ dx_1\, dx_2\le\ C_4(\ep(m)+A)\quad,
\ee
where $C_4$ is a universal constant. We set $\la_-:=\inf\{\la(x)\ ;\ x\in B_{2r}(0)\setminus B_{r}(0)\}$\,, $\la_+:=\max\{\la(x)\ ;\ x\in B_{2r}(0)\setminus B_{r}(0)\}$\,,
$\nu_-:=\inf\{\nu(x)\ ;\ x\in B_{2r}(0)\setminus B_{r}(0)\}$\,, and $\nu_+:=\max\{\nu(x)\ ;\ x\in B_{2r}(0)\setminus B_{r}(0)\}$ .
One has the estimates
\be
\label{A.a16}
\frac{e^{\la_+}}{e^{\la_-}}\le\frac{e^{\nu_++\|\mu\|_\infty}}{e^{\nu_--\|\mu\|_\infty}}\le e^{2\|\mu\|_\infty}\ e^{\nu_+-\nu_-}\le e^{2C\, \ep(m)}\ e^{2C_4 (\ep(m)+A)}\quad,
\ee
where we have used respectively (\ref{A.a11}) and (\ref{A.a15}). 

\medskip

We have thus proved (\ref{A.a3}) for $\al=1/2$. Let now $\al=2^{-j}$, for some $j\in\N^*$. The idea is to apply (\ref{A.a3}) successively $j$ times, with $\al=1/2$ and with $2^j r$ in place of $r$.  For $i=1,\ldots, j$, we let $\la^i_-:=\inf\{\la(x)\ ;\ x\in B_{2^ir}(0)\setminus B_{2^{i-1}r}(0)\}$ and
$\la^i_+:=\max\{\la(x)\ ;\ x\in B_{2^ir}(0)\setminus B_{2^{i-1}r}(0)\}$.
Using the fact that
\[
\la^i_-\le\la_+^{i\pm1}\quad,
\]
gives
\[
\frac{e^{\la^i_+}}{e^{\la^{i\pm1}_-}}\le\frac{e^{\la^i_+}}{e^{\la^{i}_-}}\,\frac{e^{\la^i_-}}{e^{\la^{i\pm1}_+}}\,\frac{e^{\la^{i\pm1}_+}}{e^{\la^{i\pm1}_-}}\le \frac{e^{\la^i_+}}{e^{\la^{i}_-}}\ 
\frac{e^{\la^{i\pm1}_+}}{e^{\la^{i\pm1}_-}}\quad.
\]
The latter enables an iteration, and eventually yields that (\ref{A.a3}) holds with $\al=2^{-j}$, for an arbitrary positive integer $j$ satisfying $4R^{-2}<2^{-j}$, thereby concluding the proof of Lemma~\ref{lm-A.0}.\hfill $\Box$

\medskip

\begin{Lma}
\label{lm-Aa0}{\bf [Pointwise control of the conformal factor in the neck region]}
There exists a constant $\eta(m)>0$ with the following property. Let $0<\eta<\eta(m)$ and   $0<4r<R<+\infty$. If $\vec{\xi}$ is any (weak) conformal immersion of $\Om:=B_{R}(0)\setminus B_r(0)$ into ${\R}^m$ with $L^2$-bounded second fundamental form, and satisfying
\be
\label{AA.b0}
\|\nabla\vec{n}_{\vec{\xi}}\|_{L^{2,\infty}(\Om)}\le\sqrt{\eta}\quad,
\ee
then there exist $\al\in (1/2,1)$ and $A\in{\R}$ (depending on $R$, $r$, $m$, and $\vec{\xi}$) such that
\be
\label{A.b01}
\|\la(x)-d\,\log|x|-A\|_{L^\infty(B_{\al R}(0)\setminus B_{\al^{-1} r}(0))}\le C_m\ \lf(\|\nabla \la\|_{L^{2,\infty}(\Om)}+\int_\Om|\nabla\vec{n}_{\vec{\xi}}|^2\, dx\rg)\quad,
\ee
where $d$ satisfies
\be
\label{AA.b01}
\begin{array}{l}
\ds\lf|2\pi\, d- \int_{\p B_r}\p_r\la\ dl_{\p B_r}\rg|\\[5mm]
\ds\quad\le C_m\ \lf[\int_{B_{2\, r}(0)\setminus B_r(0)}|\nabla\vec{n}_{\vec{\xi}}|^2\, dx+\frac{1}{\log R/r}\ \lf(\|\nabla \la\|_{L^{2,\infty}(\Om)}+\int_\Om|\nabla\vec{n}_{\vec{\xi}}|^2\, dx\rg)\rg]\quad,
\end{array}
\ee
and $C_m$ depends only on m, and $\la$ is as in lemma~\ref{lm-A.0}.
\hfill $\Box$
\end{Lma}
{\bf Proof of lemma~\ref{lm-Aa0}.}
It is convenient to introduce the constant
\[
\La:=\|\nabla \la\|_{L^{2,\infty}(\Om)}+\int_\Om|\nabla\vec{n}_{\vec{\xi}}|^2\, dx\quad.
\]
We choose $\eta(m)$ smaller than the constant $\eta(m)$ given in lemma~\ref{lm-AA.1}, so that there exists a framing $$\{\vec{e}_{1},\vec{e}_{2}\}\in (W^{1,2}(\Om,S^{m-1}))^2$$ with
\be
\label{A.bc2}
\star\,\vec{n}_{\vec{\xi}}=\vec{e}_{1}\wedge\vec{e}_{2}\quad\mbox{ and }\quad\sum_{i=1}^2\int_{\Om}|\nabla\vec{e}_{i}|^2\, dx\le C\ \int_{\Om}|\nabla\vec{n}_{\vec{\xi}}|^2\, dx\quad.
\ee
The conformal parameter $\la$ satisfies
\be
\label{A.bc3}
\Delta\la=\nabla^\perp\vec{e}_{1}\cdot\nabla\vec{e}_{2}\quad\quad\mbox{ in }\Om\quad.
\ee
Let $\mu$ be the solution of
\be
\label{A.b4}
\lf\{
\begin{array}{l}
\ds\Delta\mu=\nabla^\perp\vec{e}_{1}\cdot\nabla\vec{e}_{2}\quad\quad\mbox{ in }\Om\\[5mm]
\ds\mu=0\quad\quad\quad\mbox{ on }\p\Om\quad.
\end{array}
\rg.
\ee
As in the proof of lemma~\ref{lm-A.0}, Wente estimates give 
\be
\label{A.b5}
\|\mu\|_{L^\infty(\Om)}\le (2\pi)^{-1}\, \|\nabla\vec{e}_{1}\|_{L^2}\, \|\nabla\vec{e}_{2}\|_{L^2}\le C\, \La
\ee
and 
\be
\label{A.b6}
\|\nabla\mu\|_{L(\Om)}\le \sqrt{3/32\pi}\, \|\nabla\vec{e}_{1}\|_{L^2}\, \|\nabla\vec{e}_{2}\|_{L^2}\le C\, \La\quad.
\ee
Let $\nu:=\la-\mu$ be the harmonic function in $\Om$ equal to $\la$ on the boundary $\p\Om$. Analogously to (\ref{A.ba3}) in the proof of  lemma~\ref{lm-A.a00}, there exists $\al\in(1/2,1)$ with
\be
\label{A.ba7}
\begin{array}{l}
\ds\int_{\p (B_{\al R}\setminus B_{\al^{-1} r})}|\nabla \nu|\le r\ \|\nabla \nu\|_{L^1(B_{2\,r}\setminus B_{r})}
+ 2\,R\ \|\nabla \nu\|_{L^1(B_{R}\setminus B_{R/2})}\\[5mm]
\ds\quad\quad\quad\le C_4\, \|\nabla \nu\|_{L^{2,\infty}(\Om)}\le C_4\, \lf[\|\nabla\la\|_{L^{2,\infty}(\Om)}+\|\nabla\mu\|_{L(\Om)}\rg]
\le C_4\,(1+C)\, \La\quad,
\end{array}
\ee
where $C_4$ is a universal constant. From (\ref{A.ba7})
we deduce the existence of two constants $a$ and $b$ depending on $\vec{\xi}$ such that
\be
\label{A.ba8}
\forall\: x\in \p B_{\al^{-1}r}(0)\quad\quad\quad-C_m\,\La\le \nu(x)-a\le C_m\,\La\quad,
\ee
and
\be
\label{A.ba9}
\forall\: x\in \p B_{\al R}(0)\quad\quad\quad-C_m\,\La\le \nu(x)-b\le C_m\,\La\quad,
\ee
where $C_m>0$ only depends on $m$.  Since the functions
\[
\tau^\pm(x):=\frac{b-a}{\log R/r}\log\frac{|x|}{r}+a\pm C_m\,\La
\]
are harmonic, and since (\ref{A.ba8}) and (\ref{A.ba9}) imply
\[
\tau^-\le\nu\le\tau^+\quad\quad\quad\mbox{on }\p( B_{\al R}\setminus B_{\al^{-1} r})\quad,
\]
the maximum principle yields
\be
\label{A.ba9a}
\forall\: x\in B_{\al R}\setminus B_{\al^{-1} r}\quad\quad\quad-C_m\,\La\le\nu(x)-\lf[\frac{b-a}{\log R/r}\log\frac{|x|}{r}+a\rg]\le C_m\,\La\quad.
\ee
Let  $\rho\in (r,R)$, and define $2\pi\,d:=\int_{\p B_\rho}\p_\rho\nu$. Integrating by parts gives
\[
\begin{array}{l}
\ds 0=\int_{B_{\al R}\setminus B_{\al^{-1} r}}\Delta\nu\ \log|x|\, dx\\[5mm]
\ds\quad=2\pi\,d\ \log\frac{R}{r}+\int_{0}^{2\pi}\nu(\al^{-1}\,r,\theta)\ d\theta-\int_{0}^{2\pi}\nu(\al\,R,\theta)\ d\theta\quad.
\end{array}
\]
Hence
\be
\label{A.ba10}
\lf|d-\frac{b-a}{\log R/r}\rg|\le \frac{C_m\ \La}{\log{R/r}}\quad. 
\ee
Combining (\ref{A.ba9a}) and (\ref{A.ba10}) yields
\be
\label{A.ba11}
\forall\: x\in B_{\al R}\setminus B_{\al^{-1} r}\quad\quad\quad-C_m\,\La\le\nu(x)-\frac{1}{2\pi}\lf[\int_{\p B_{r}(0)}\p_r\nu\rg]\, \log\frac{|x|}{r}-a\le C_m\,\La\quad.
\ee

\medskip

We next estimate the difference $\int_{\p B_{r}(0)}\p_r\mu=\int_{\p B_{r}(0)}\p_r(\nu-\la)$. There holds
\be
\label{A.ba12}
\begin{array}{l}
\ds 0=\int_{\Om}\mu\ \Delta\log({|x|}/{R})\, dx=-\int_{\Om}\nabla\mu\cdot\nabla\log({|x|}/{R})\, dx\\[5mm]
\ds=-\log({r}/{R})\ \int_{\p B_r}\p_r\mu+
\int_\Om(\nabla^\perp\vec{e}_1\cdot\nabla\vec{e}_2)\ \log({|x|}/{R})\, dx\quad.
\end{array}
\ee
As seen in the proof of lemma~\ref{lm-AA.1}, the construction of the moving frame $\{\vec{e}_1,\vec{e}_2\}$ requires to construct an extension $\hat{\vec{n}}$ of $\vec{n}_{\vec{\xi}}$
into $B_r(0)$. Since $\{\vec{e}_1,\vec{e}_2\}$ is also a moving frame for $\hat{\vec{n}}$ in $B_R(0)$, lemma~\ref{lm-A.2} shows that
\[
\sum_{i=1}^2\int_{B_{R}(0)}|\nabla\vec{e}_{i}|^2\, dx\le C\,\int_{\Om}|\nabla\vec{n}_{\vec{\xi}}|^2\, dx\le \La\quad.
\]
Using this extension, we may also recast (\ref{A.ba12}) in form
\be
\label{A.ba13}
\begin{array}{l}
\ds\log({R}/{r})\ \int_{\p B_r}\p_r\mu=\int_{B_R}(\nabla^\perp\vec{e}_1\cdot\nabla\vec{e}_2)\ \log({|x|}/{R})\, dx\\[5mm]
\ds\quad\quad\quad\quad-\int_{B_r}(\nabla^\perp\vec{e}_1\cdot\nabla\vec{e}_2)\ \log({|x|}/{R})\, dx\quad.
\end{array}
\ee
Let $\Psi$ be the solution of
\[
\lf\{
\begin{array}{l}
\ds\Delta\Psi=\nabla^\perp\vec{e}_1\cdot\nabla\vec{e}_2\quad\quad\mbox{ in }B_R(0)\\[5mm]
\ds\Psi=0\quad\quad\quad\mbox{ on }\p B_R(0)\quad.
\end{array}
\rg.
\]
As in the proof of lemma~\ref{lm-A.0}, Wente estimates give
\be
\label{A.ba14}
\|\nabla\Psi\|_{L^{2,1}(B_R(0))}\le C_0\ \|\nabla\vec{e}_1\|_{L(B_R(0))} \|\nabla\vec{e}_2\|_{L(B_R(0))}\le C_m\ \La\quad,
\ee
for some universal constant $C_0$ independent of $R$. Hence,
\be
\label{A.ba15} 
\begin{array}{l}
\ds\lf|\int_{B_R}(\nabla^\perp\vec{e}_1\cdot\nabla\vec{e}_2)\ \log({|x|}/{R})\, dx\rg|=\lf|\int_{B_R}\nabla\Psi\ \nabla\log({|x|}/{R})\, dx\rg|\\[5mm]
\ds\quad\quad\le\|\nabla\Psi\|_{L^{2,1}(B_R(0))}\ \|\nabla\log({|x|}/{R})\|_{L^{2,\infty}(B_R(0))}\le C_m\ \La\quad.
\end{array}
\ee
On the other hand, lemma~\ref{lm-AA.1} and lemma~\ref{lm-A.2} yield a frame $$\{\vec{f}_{1},\vec{f}_{2}\}\in (W^{1,2}(B_{r}(0),S^{m-1}))^2$$ satisfying
$$
\star\,\hat{\vec{n}}=\vec{f}_{1}\wedge\vec{f}_{2}
$$
and
\be
\label{A.b2}
\sum_{i=1}^2\int_{B_{r}(0)} |\nabla\vec{f}_{i}|^2\, dx\le C\ \int_{B_{r}(0)}|\nabla\hat{\vec{n}}|^2\, dx\le C_m\ \int_{B_{2r}(0)\setminus B_r(0)}|\nabla\vec{n}_{\vec{\xi}}|^2\, dx\quad.
\ee
Clearly, there exists $\theta\in W^{1,2}(B_r,{\R})$ such that
\[
e^{i\theta} (\vec{e}_1+i\vec{e}_2)=\vec{f}_1+i\vec{f}_2\quad,
\]
and thus
\[
\vec{e}_1\cdot\nabla\vec{e}_2\,+\,\nabla\theta=\vec{f}_1\cdot\nabla\vec{f}_2\quad\quad\mbox{ in }B_{r}(0)\quad,
\]
Whence there holds
\be
\label{A.b3}
\nabla^\perp\vec{e}_{1}\cdot\nabla\vec{e}_{2}=\nabla^\perp\vec{f}_{1}\cdot\nabla\vec{f}_{2}\quad\quad\mbox{ in }B_{r}(0)\quad.
\ee
Let $\psi$ be the solution of 
\[
\lf\{
\begin{array}{l}
\ds\Delta\psi=\nabla^\perp\vec{f}_1\cdot\nabla\vec{f}_2\,-\,|B_r|^{-1}\,\int_{B_r} \nabla^\perp\vec{f}_1\cdot\nabla\vec{f}_2\quad\quad\mbox{ in }B_r(0)\\[5mm]
\ds\p_r\psi=0\quad\quad\quad\mbox{ on }\p B_r(0)\quad.
\end{array}
\rg.
\]
Calling again upon Wente estimates (this time with Neuman boundary data) gives the following control of the $L^{2,1}$-norm of $\nabla\psi$\,:
\be
\label{A.ba16}
\|\nabla\psi\|_{L^{2,1}(B_R(0))}\le C_0\ \|\nabla\vec{f}_1\|_{L(B_R(0))} \|\nabla\vec{f}_2\|_{L(B_R(0))}\le C_m\ \int_{B_{2r}(0)\setminus B_r(0)}|\nabla\vec{n}_{\vec{\xi}}|^2\, dx\quad.
\ee
Consequently, we derive
\be
\label{A.ba17} 
\begin{array}{l}
\ds\lf|\int_{B_r}(\nabla^\perp\vec{e}_1\cdot\nabla\vec{e}_2)\ \log({|x|}/{R})\, dx\rg|\le\log({R}/{r})\lf|\int_{B_r}\nabla^\perp\vec{f}_1\cdot\nabla\vec{f}_2\rg|+\lf|\int_{B_r}\nabla\psi\ \nabla\log({|x|}/{R})\, dx\rg|\\[5mm]
\ds\quad\quad\le C_m\ \log({R}/{r})\int_{B_{2r}(0)\setminus B_r(0)}|\nabla\vec{n}_{\vec{\xi}}|^2\, dx+
\|\nabla\psi\|_{L^{2,1}(B_r(0))}\ \|\nabla\log({|x|}/{R})\|_{L^{2,\infty}(B_R(0))}\\[5mm]
\ds\quad\quad\le C_m\ \log({R}/{r})\int_{B_{2r}(0)\setminus B_r(0)}|\nabla\vec{n}_{\vec{\xi}}|^2\, dx\quad.
\end{array}
\ee
Finally, combining (\ref{A.ba13}), (\ref{A.ba15}), and (\ref{A.ba17}) gives the announced
\be
\label{A.ba18}
\log({R}/{r})\ \lf|\int_{\p B_r}\p_r\mu\rg|\le C_m\ \lf[\La+\log({R}/{r})\ \int_{B_{2r}(0)\setminus B_r(0)}|\nabla\vec{n}_{\vec{\xi}}|^2\, dx\rg]\quad.
\ee
\hfill $\Box$

\section{Uniform $L^{2,1}$-control of the mean curvature vector of a Willmore immersion in neck regions   }
\reset

The goal of this section is to prove the following result.

\begin{Lm}
\label{lm-VII.1}
Let $m\ge 3$. There exists a constant $\ep(m)>0$ with the following property. Let $0<8r<R$, and let $\vec{\xi}$ be any conformal Willmore immersion of $B_R(0)$ into ${\R}^m$ such that
\be
\label{VII.1}
\sup_{r<s<R/2} \int_{B_{2s}(0)\setminus B_s(0)}|\nabla\vec{n}_{\vec{\xi}}|^2\, dx\le \ep(m)\quad.
\ee
We set
\[
\La\,:=\,\|\nabla \la\|_{L^{2,\infty}(B_R\setminus B_r)}+\int_{B_R(0)}|\nabla\vec{n}_{\vec{\xi}}|^2\, dx\quad,
\]
where $\la$ denotes as before the conformal parameter of $\vec{\xi}$. Then there holds 
\be
\label{VII.2}
\|e^{\la}\vec{H}_{\vec{\xi}}\|_{L^{2,1}(B_{R/2}\setminus B_{2r})}\le C(m,\La)\quad,
\ee
where $\vec{H}_{\vec{\xi}}$ is the mean curvature vector of the immersion $\vec{\xi}$, and $C(m,\La)$ is a positive constant depending\footnote{the constant $C(m,\La)$ is independent of the conformal type of the annulus.} only on $m$ and $\La$.\hfill $\Box$
\end{Lm}
{\bf Proof of lemma~\ref{lm-VII.1}.}
According to the work in \cite{Ri2}, an immersion $\vec{\xi}$ is conformal and Willmore on $B_R(0)$ if and only if there exists an ${\R}^m$-valued map $\vec{L}$ (uniquely defined up to an additive constant vector) satisfying
\be
\label{VII.3}
\nabla\vec{L}=2\nabla^{\perp}\vec{H}-3\nabla^\perp(\pi_{\vec{n}})(\vec{H})+\star\,(\nabla\vec{n}_{\vec{\xi}}\wedge\vec{H})\quad,
\ee
where $\star$ is the canonical Hodge operation on multivectors of ${\R}^m$ and $\pi_{\vec{n}}$ is the projection operator onto the normal space to $\vec{\xi}_\ast T\Sigma$. It is given by two consecutive applications of the contraction operator\footnote{The contraction operator $\res$ is the linear operation defined as follows. For every $\al$, $\beta$, and $\gamma$, respectively $p$, $q$, and $(p-q)$ vectors in ${\R}^m$, 
one has
\[
<\al\res\beta,\gamma>=<\al,\beta\wedge\gamma>\quad.
\]}
$\res$ with $\vec{n}$ (see (VI.87) in \cite{Ri1})\,:
\[
\forall\: \vec{w}\in{\R}^m\quad\quad\pi_{\vec{n}}(\vec{w}):=(-1)^{m-1}\vec{n}\res(\vec{n}\res\vec{w})\quad.
\]
Hence the gradient of $\pi_{\vec{n}}$ in (\ref{VII.3}) is to be understood as
\be
\label{VII.4}
\forall\: \vec{w}\in{\R}^m\quad\quad\nabla(\pi_{\vec{n}})(\vec{w}):=(-1)^{m-1}(\nabla\vec{n})\res(\vec{n}\res\vec{w})+(-1)^{m-1}\vec{n}\res((\nabla\vec{n})\res\vec{w})\quad.
\ee

The first condition to impose on the constant $\ep(m)$ is that it should be smaller than the one in lemma~\ref{lm-A.0}, so that the {\it uniform Harnack estimate} (\ref{A.a3}) holds on the annulus $B_{\al R}\setminus B_{\al^{-1} r}$ for every $0<\al<1$. If, in addition, we ensure that $\ep(m)$ be smaller than the threshold constant appearing the $\epsilon$-regularity theorem I.5 in \cite{Ri2}, then the uniform Harnack estimate enables\footnote{indeed, the proof of theorem I.5 in \cite{Ri2} requires that the conformal parameter $\la$ be ``roughly" constant on balls around the origin. We are instead working on annuli, so the estimate (\ref{A.a3}) is exactly what is needed.} us to conclude that for any $x\in B_{4R/5}(0)\setminus B_{5\,r/4}(0)$ there holds
\be
\label{VII.5}
|\nabla\vec{n}_{\vec{\xi}}(x)|^2\le C\,|x|^{-2}\int_{B_{2|x|}(0)\setminus B_{|x|/2}(0)}|\nabla\vec{n}_{\vec{\xi}}|^2\, dx\le C\,\ep(m)\,x^{-2}\quad,
\ee
where $C$ is some constant independent of the data of the problem. 
This implies in particular that
\be
\label{VII.5a}
\|\nabla\vec{n}_{\vec{\xi}}\|_{L^{2,\infty}( B_{4R/5}(0)\setminus B_{5\,r/4}(0))}\le \sqrt{\frac{15\pi}{16}\,C\,\ep(m)}=C'\,\sqrt{\ep(m)}\quad.
\ee
Finally, we ask that $\ep(m)>0$ satisfy $C'\sqrt{\ep(m)}\le \sqrt{\eta(m)}$, where $\eta(m)$ is the constant given in lemma~\ref{lm-Aa0}. Accordingly, there exist a constant $A$ (which depends a-priori on all parameters $r$, $R$, $\La$, and $\vec{\xi}$) and a constant $d$ such that
\be
\label{VII.6a}
|d|\le C_m\ \La\quad\,
\ee
and
\be
\label{VII.6b}
\|\la(x)-d\,\log|x|-A\|_{L^\infty(B_{R/2}\setminus B_{2r})}\le C_m\,\La\quad.
\ee
We introduce for $0<s<R/2$ the function 
$$\delta(s):=\lf(s^{-2}\int_{B_{2s}(0)\setminus B_{s/2}(0)}|\nabla\vec{n}_{\vec{\xi}}|^2\,dx\rg)^{1/2}\quad.$$
Clearly,
\be
\label{triv}
s\,\delta(s)\le \sqrt{\La}\qquad\quad\forall\: s\in(0,R/2)\quad.
\ee
Inequality (\ref{VII.5}) gives
\be
\label{VII.6c}
e^{\la(x)}\,|\vec{H}_{\vec{\xi}}(x)|\le |\nabla\vec{n}_{\vec{\xi}}(x)|\le C\,\delta(|x|)\quad\qquad\forall\: x\in B_{R/2}(0)\setminus B_{2r}(0)\quad.
\ee
Clearly, since $8r<R$, there holds
\be
\label{VII.6}
\int_{r}^{R/2}s\,\delta^2(s)\,ds\le \log(4)\,\int_{B_{R}(0)\setminus B_{r/2}(0)}|\nabla\vec{n}_{\vec{\xi}}|^2\, dx\le 2\,\La\quad.
\ee
The identity (\ref{VII.3}) shows that on any dyadic annulus $B_{2s}\setminus B_s$ included in $B_{R}(0)\setminus B_r(0)$, the mean curvature vector satisfies
\[
\Delta\vec{H}=\frac{1}{2}\,div\lf(3\nabla(\pi_{\vec{n}})(\vec{H})+\star(\nabla^\perp\vec{n}_{\vec{\xi}}\wedge\vec{H})\rg)\quad.
\]
From (\ref{VII.6c}), a standard argument from the theory of second-order uniformly elliptic equations in divergence form then yields the existence of a constant $C_1>0$ such that
\be
\label{VII.8}
e^{\la(x)}\,|\nabla\vec{H}_{\vec{\xi}}(x)|\le C_1\, \frac{\delta(|x|)}{|x|}\qquad\quad \forall\: x\in B_{R/2}(0)\setminus B_{2r}(0)\quad.
\ee
For any $t\in (r,R)$, we denote
\[
\vec{L}_t:=\frac{1}{|\p B_t(0)|}\int_{\p B_t(0)}\vec{L}\quad,
\]
where the function $\bL$ is as (\ref{VII.6c}). Using (\ref{VII.3}), (\ref{VII.6b}), and (\ref{VII.8}), there holds
\be
\label{VII.9}
|\vec{L}(x)-\vec{L}_{|x|}|\le\int_{\p B_{|x|}}|\nabla\vec{L}|\le 2\pi\,C_{1}\,e^{2C_m\La}\, e^{-\la(x)}\,\delta(|x|)\quad\qquad \forall\:  x\in B_{R/2}(0)\setminus B_{2r}(0)\quad.
\ee
This implies in particular that
\be
\label{VII.10}
\int_{B_{R/2}(0)\setminus B_{2r}(0)} e^{2\la(x)}\, |\vec{L}(x)-\vec{L}_{|x|}|^2\, dx\le\, C_2\,\La\quad,
\ee
for some constant $C_2$. Note that we have used (\ref{VII.6}). \\
We have next
\be
\label{VII.11}
\frac{d \vec{L}_t}{dt}=\frac{1}{2\pi}\int_0^{2\pi}\ \frac{\p \vec{L}}{\p t}(t,\theta)\ d\theta=\frac{3}{2\pi}\int_0^{2\pi}\frac{1}{r}\frac{\p \pi_{\vec{n}}}{\p\theta}(\vec{H})\ d\theta+
\frac{1}{2\pi}\int_0^{2\pi}\star\lf(\frac{\p \vec{n}_{\vec{\xi}}}{\p \nu}\wedge\vec{H}\rg)\ d\theta\quad.
\ee
Setting $a(t):=|\vec{L}_t|$, and putting (\ref{VII.6c}) into the latter gives the estimate
\be
\label{VII.12}
|\dot{a}(t)|=\lf|\frac{d |\vec{L}_t|}{dt}\rg|= \lf|\frac{d \vec{L}_t}{dt}\rg|\le C^2\, e^{-\la}\, \delta^2(t)\quad,
\ee
Whence, using (\ref{VII.6}), we reach
\be
\label{VII.13}
\int_{2r}^{R/2} s\, e^\la\, |\dot{a}(s)|\, ds\le 2\,C^2 \La\quad.
\ee
The Harnack estimate (\ref{VII.6b}) yields
\be
\label{VII.14}
e^{A-C_m\La}\,  |x|^d\, \le e^{\la(x)}\le e^{A+C_m\La}\,  |x|^d\quad\qquad\forall\: x\in B_{R/2}\setminus B_{2r}(0) \quad,
\ee
thereby showing that (\ref{VII.13}) may be recast in the form
\be
\label{VII.15}
\int_{2r}^{R/2} s^{1+d}\, |\dot{a}(s)|\, ds\le 2\,C^2\, e^{C_m\,\La-A}\,\La \quad.
\ee
An elementary integration by parts gives for any $r\le \tau<T\le R$\,:
\[
\int_\tau^Ts^{1+d}\, \dot{a}(s)\, ds=T^{1+d}\, a(T)-\tau^{1+d}\, a(\tau)-(1+d)\int_{\tau}^Ts^d\, a(s)\, ds\quad.
\]
Hence, since $a\ge 0$, we have
\be
\label{VII.16}
\lf\{
\begin{array}{l}
\ds\forall\: \ d\le -1 \quad\quad \forall\:\  2r<t<R/2\quad\quad t^{1+d}\, a(t)\le (2r)^{1+d}\, a(2r)+\int_{2r}^{R/2} s^{1+d}\ |\dot{a}(s)|\, ds\\[5mm]
\ds\forall\:\  d\ge -1 \quad\quad \forall\:\  2r<t<R/2\quad\quad t^{1+d}\, a(t)\le (R/2)^{1+d}\, a(R/2)+\int_{2r}^{R/2} s^{1+d}\ |\dot{a}(s)|\, ds\quad.
\end{array}
\rg.
\ee
Recall that we still have the freedom to adjust the vector $\vec{L}$ by an additive constant vector. It will be convenient to us to choose $\vec{L}$ in such a way that
\[
\mbox{ If } \ d\le -1 \mbox{ we take }\int_{\p B_{2r}}\vec{L}=0\quad\mbox{ whereas if } \ d\ge -1 \mbox{ we take }\int_{\p B_{R/2}}\vec{L}=0\quad.
\]
This particular choice implies
\be
\label{VII.17}
t\, e^\la\, |\vec{L}_t|\equiv t\, e^\la\, a(t)\le 2\,C^2\, e^{2\, C_m\,\La}\, \La\quad\qquad\forall\:\  2r<t<R/2\quad,
\ee
where we have used (\ref{VII.16}), (\ref{VII.14}) and (\ref{VII.15}). With the help of an elementary computation, the latter then implies that
\be
\label{VII.19}
\|e^{\la(x)}|\vec{L}_{|x|}|\|_{L^{2,\infty}(B_{R/2}\setminus B_{2r})}\le C_0(m,\La)\quad,
\ee
for some positive constant $C_0(m,\La)$. 
Combining this inequality with (\ref{VII.10}) gives
\be
\label{VII.20}
\|e^{\la(x)}\vec{L}(x)\|_{L^{2,\infty}(B_{R/2}\setminus B_{2r})}\le C(m,\La)\quad,
\ee
where $C(m,\La)$ depends only on $m$ and $\La$. Moreover, the estimates (\ref{VII.9}) and (\ref{VII.17}) altogether yield
\be
\label{VII.20a}
e^{\la(x)}|\vec{L}(x)|\le C(\La,C_m)\,\Big[1+|x|\delta(x)\Big]|x|^{-1}\le C'(m,\La)\, |x|^{-1}\quad\qquad \forall\: x\in B_{R/2}(0)\setminus B_{2r}(0)\quad.
\ee
Note that we have also used (\ref{triv}).

\medskip

As is done in \cite{Ri2} (see also section VI.7.2 of \cite{Ri1}), we introduce on $B_R(0)$ the smooth function $S$ and the smooth $\bigwedge^{2}(\R^m)$-valued map $\vec{R}$ via
\be
\label{VI.172a}
\lf\{
\begin{array}{l}
\ds \nabla S=\vec{L}\cdot\nabla\vec{\xi}\\[2mm]
\ds\nabla\vec{R}=\vec{L}\wedge\nabla\vec{\xi}+2\ (\star(\vec{n}_{\xi}\res\vec{H}))\res\nabla\vec{\xi}\quad.
\end{array}
\rg.
\ee
Since $|\nabla\vec{\xi}|^2=2\,e^{2\la}$, the estimates (\ref{VII.5a}), (\ref{VII.20}) give the bound
\be
\label{VI.21}
\|\nabla S\|_{L^{2,\infty}(B_{R/2}\setminus B_{2r})}+\|\nabla \vec{R}\|_{L^{2,\infty}(B_{R/2}\setminus B_{2r})}\le C(m,\La)\quad.
\ee
Furthermore, using (\ref{VII.6c}), (\ref{VII.20a}), and (\ref{triv}), we obtain the pointwise bound
\be
\label{VII.21a}
|\nabla S(x)|+|\nabla\vec{R}(x)|\le C(m,\La)\,|x|^{-1}\quad\qquad\forall\: x\in B_{R/2}\setminus B_{2r}\quad.
\ee

\medskip

\noindent One verifies (cf. \cite{Ri2} and theorem VI.15 of \cite{Ri1}) that the following equations\footnote{The linear operator $\bullet$ is the contraction 
which to a pair of $p$ and $q$-vectors of ${\R}^m$ assigns
a $(p+q-2)$-vector of ${\R}^m$ in such a way that
\[
\forall\:\,\vec{a}\in\wedge^p{\R}^m\quad\ \forall\:\,\vec{b}\in\wedge^1{\R}^m\quad\vec{a}\bullet\vec{b}:=\vec{a}\res\vec{b}
\]
and
\[
\begin{array}{l}
\ds\forall\:\,\vec{a}\in\wedge^p{\R}^m\quad\ \forall\:\,\vec{b}\in\wedge^r{\R}^m\quad\ \forall\:\,\vec{c}\in\wedge^s{\R}^m\quad\\[5mm]
\ds\vec{a}\bullet(\vec{b}\wedge\vec{c}):=(\vec{a}\bullet\vec{b})\wedge\vec{c}+(-1)^{r\,s}(\vec{a}\bullet\vec{c})\wedge\vec{b}\quad.
\end{array}
\]
}
hold
\be
\label{VI.173}
\lf\{
\begin{array}{rcl}
\ds\nabla S&=&\ds-(\star\,\vec{n})\cdot\nabla^\perp\vec{R}\\[2mm]
\ds\nabla\vec{R}&=&\ds(-1)^m\ \star(\vec{n}\bullet\nabla^\perp\vec{R})+ (\star\,\vec{n})\nabla^\perp S\quad.
\end{array}
\rg.
\ee
For any $t\in(2r, R/2)$, we let
\[
S_t:=\frac{1}{|\p B_t(0)|}\int_{\p B_t(0)} S\quad\quad\mbox{,}\quad\quad \vec{R}_t:=\frac{1}{|\p B_t(0)|}\int_{\p B_t(0)} \vec{R}\quad\quad\mbox{and}\quad\quad\vec{n}_t:=\frac{1}{|\p B_t(0)|}\int_{\p B_t(0)} \vec{n}\quad.
\]
The equations (\ref{VI.173}) then yield
\be
\label{VII.22}
\lf\{
\begin{array}{l}
\ds \frac{d S_t}{dt}=\frac{1}{2\pi}\int_0^{2\pi}\frac{\p S}{\p t}(t,\theta)\,d\theta=\frac{1}{2\pi}\int_0^{2\pi} \star\,(\vec{n}-\vec{n}_t)\cdot\frac{1}{t} \frac{\p \vec{R}}{\p\theta}\ d\theta\\[6mm]
\ds \frac{d \vec{R}_t}{dt}=\frac{1}{2\pi}\int_0^{2\pi}\frac{\p \vec{R}}{\p t}(t,\theta)\,d\theta\\[5mm]
\ds\quad\quad=(-1)^{m-1}\frac{1}{2\pi}\int_0^{2\pi}(\vec{n}-\vec{n}_t)\bullet\frac{1}{t} \frac{\p \vec{R}}{\p\theta}\ d\theta
-\frac{1}{2\pi}\int_0^{2\pi}\star\,(\vec{n}-\vec{n}_t)\,\frac{1}{t} \frac{\p S}{\p\theta}\, d\theta
\end{array}
\rg.
\ee
Note that (\ref{VII.6c}) gives
\be
\label{VII.23}
|\vec{n}(x)-\vec{n}_{|x|}|\le C\, |x|\, \delta(|x|)\quad\qquad \forall\: x\in B_{R/2}(0)\setminus B_{2r}(0)\quad.
\ee
Thus, (\ref{VII.21a}), (\ref{VII.22}), and (\ref{VII.23}) altogether yield
\be
\label{VII.24}
\lf|\frac{d S_t}{dt}\rg| +\lf|\frac{d \vec{R}_t}{dt}\rg|\le\, C(m,\La)\ \delta(t)\quad.
\ee
Hence, with the help of (\ref{VII.6}), we deduce
\be
\label{VII.25}
\int_{2r}^{R/2}\lf[\lf|\frac{d S_t}{dt}\rg|^2 +\lf|\frac{d \vec{R}_t}{dt}\rg|^2\rg]\,t\,dt\le C(m,\La)\quad.
\ee
Our study requires the following result, proved in \cite{LaRi}. 
\begin{Lm}
\label{lm-VII.20}\cite{LaRi}
Let $a$ and $b$ be two functions on $B_1(0)$ such that $\nabla a\in L^{2,\infty}$ and $\nabla b\in L^2$. Let  $0<\ep<1/4$ and $\phi$ satisfy
\[
-\Delta\phi=\p_{x_1}a\,\p_{x_2} b-\p_{x_1}b\,\p_{x_2} a\quad\quad\mbox{ in }B_1(0)\setminus B_\ep(0)\quad,
\]
For $\ep\le r\le1$, we set $\,\phi_0(r):=(2\pi\, r)^{-1}\int_{\p B_r(0)}\phi$, and we assume that
\be
\label{x1}
\int_{\ep}^1|\dot{\phi_0}|^2\,r\,dr<+\infty\quad.
\ee
Then $\nabla\phi\in L^{2}(B_{1/2}(0)\setminus B_{2\ep}(0))$ and  there exists a positive constant $C$ independent of $\ep$ and $\phi$ such that 
\be
\label{x2}
\begin{array}{l}
\ds\| \nabla \phi \|_{L^{2}(B_{1/2}(0)\setminus B_{2\ep}(0))} \leq C\ \| \nabla a\|_{{2,\infty}}\,  \|\nabla b\|_{2} + C\, \|\nabla\phi_0\|_{L^2(B_1(0)\setminus B_\ep(0))}\\[3mm]
\ds\hspace{4cm} \quad\quad+\: C\, \|\nabla \phi\|_{L^{2,\infty}(B_1(0)\setminus B_\ep(0))}\quad.
\end{array}
\ee
\hfill $\Box$
\end{Lm}
Applying divergence to both sides of the equations (\ref{VI.173}) gives the {\it conservative conformal Willmore system}\,:
\be
\label{VII.26}
\lf\{
\begin{array}{rcl}
\ds\Delta S&=&\ds-\star\nabla\vec{n}\cdot\nabla^\perp\vec{R}\\[5mm]
\ds\Delta\vec{R}&=&\ds(-1)^m \star(\nabla\vec{n}\bullet\nabla^\perp\vec{R})+\star\nabla\vec{n}\cdot\nabla^\perp S \quad.
\end{array}
\rg.
\ee
Owing to (\ref{VII.1}), (\ref{VI.21}), and (\ref{VII.25}), we may apply lemma \ref{lm-VII.20} to the system (\ref{VII.26}), and find
\be
\label{VII.277}
\|\nabla S\|_{L^2(B_{R/4}(0)\setminus B_{4r}(0))}+\|\nabla \vec{R}\|_{L^2(B_{R/4}(0)\setminus B_{4r}(0))}\le\ C(m,\La)\quad.
\ee
Since the $L^{2}$-norm of $|x|^{-1}$ on annuli of the form $B_{as}(0)\setminus B_{bs}(0)$ is independent of $s$, we can use the pointwise estimate (\ref{VII.21a}) to get an upper bound independent of $R$ and $r$ for 
$$
\|\nabla S\|_{L^2(B_{R/2}(0)\setminus B_{R/4}(0)\cup B_{4r}(0)\setminus B_{2r}(0))}+\|\nabla \vec{R}\|_{L^2(B_{R/2}(0)\setminus B_{R/4}(0)\cup B_{4r}(0)\setminus B_{2r}(0))}\quad.
$$
Combining it to (\ref{VII.277}) then yields
\be
\label{VII.27}
\|\nabla S\|_{L^2(B_{R/2}(0)\setminus B_{2r}(0))}+\|\nabla \vec{R}\|_{L^2(B_{R/2}(0)\setminus B_{2r}(0))}\le\ C(m,\La)\quad.
\ee
This new information, along with (\ref{VII.6}), is injected into (\ref{VII.22}) so as to produce
\be
\label{VII.29}
\begin{array}{l}
\ds\int_{2r}^{R/2}\lf|\frac{d S_t}{dt}\rg| +\lf|\frac{d \vec{R}_t}{dt}\rg|\,dt\le\  C\int_{2r}^{R/2} \delta(t)\, \int_{\p B_t}|\nabla S|+|\nabla \vec{R}|\, dt\\[5mm]
\ds\quad\quad\le C\ \lf(\int_{2r}^{R/2}\delta^2(t)\, t\, dt\rg)^{1/2}\ \lf(\int_{B_{R/2}\setminus B_{2r}}|\nabla S|^2+|\nabla \vec{R}|^2\, dx\rg)^{1/2}\le C(m,\La)\quad.
\end{array}
\ee
Since the functions $S$ and $\bR$ are both defined up to an additive ``constant", we have the freedom to impose the conditions $S_{2r}=0$ and $\vec{R}_{2r}=\vec{0}$. From (\ref{VII.24}), we deduce in particular that
\be
\label{VII.30}
|S_t|+|\vec{R}_t|\le C(m,\La)\quad\qquad \forall\: t\in (2r,R/2)\quad.
\ee
Paired to the pointwise estimate (\ref{VII.21a}) on the gradient of $S$ and ${\vec{R}}$, the latter implies
\be
\label{VII.31}
\|S\|_{L^\infty(B_{R/2}(0)\setminus B_{2r}(0))}+\|\vec{R}\|_{L^\infty(B_{R/2}(0)\setminus B_{2r}(0))}\le\ C(m,\La)
\ee
We are sufficiently geared to apply the following general result, whose proof may be found in \cite{LaRi}.
\begin{Lm}
\label{lm-VII.2}\cite{LaRi}
Let $a$ and $b$ be two functions on $B_1(0)$ such that $\nabla a\in L^{2}$ and $\nabla b\in L^2$. Let  $0<\epsilon<1/4$ and $\phi$ satisfy
\[
-\Delta\phi=\p_{x_1}a\,\p_{x_2} b-\p_{x_1}b\,\p_{x_2} a\quad\quad\mbox{ in }B_1(0)\setminus B_\ep(0)\quad.
\]
Assume that
\be
\label{x4}
\|\phi\|_\infty<+\infty\quad.
\ee
Then $\nabla\phi\in L^{2,1}(B_{1/2}(0)\setminus B_{2\ep}(0))$, and there exists a positive constant $C$ independent of $\ep$ and $\phi$ such that 
\be
\label{x5}
\begin{array}{l}
\ds\| \nabla \phi \|_{L^{2,1}(B_{1/2}(0)\setminus B_{2\ep}(0))} \leq C\, \| \nabla a\|_{{2}}\,  \|\nabla b\|_{2} + C\, \|\phi\|_{L^\infty(B_1(0)\setminus B_\ep(0))}\\[3mm]
\ds\hspace{4cm}\quad\quad+ \:C\, \|\nabla \phi\|_{L^{2}(B_1(0)\setminus B_\ep(0)}\quad.
\end{array}
\ee
\hfill $\Box$
\end{Lm}
Just as we did above, owing to (\ref{VII.1}), (\ref{VI.21}), and (\ref{VII.31}), we apply lemma \ref{lm-VII.2} to the system (\ref{VII.26}), and find
\be
\label{VII.327}
\|\nabla S\|_{L^{2,1}(B_{R/4}(0)\setminus B_{4r}(0))}+\|\nabla \vec{R}\|_{L^{2,1}(B_{R/4}(0)\setminus B_{4r}(0))}\le\ C(m,\La)\quad.
\ee
Since the $L^{2,1}$-norm of $|x|^{-1}$ on annuli of the form $B_{as}(0)\setminus B_{bs}(0)$ is independent of $s$, we use (\ref{VII.21a}) to get an upper bound independent of $R$ and $r$ for 
$$
\|\nabla S\|_{L^{2,1}(B_{R/2}(0)\setminus B_{R/4}(0)\cup B_{4r}(0)\setminus B_{2r}(0))}+\|\nabla \vec{R}\|_{L^{2,1}(B_{R/2}(0)\setminus B_{R/4}(0)\cup B_{4r}(0)\setminus B_{2r}(0))}\quad.
$$
Combining it to (\ref{VII.327}) then yields
\be
\label{VII.32}
\|\nabla S\|_{L^{2,1}(B_{R/2}(0)\setminus B_{2r}(0))}+\|\nabla \vec{R}\|_{L^{2,1}(B_{R/2}(0)\setminus B_{2r}(0))}\le\ C(m,\La)\quad.
\ee
It is shown in \cite{Ri2} that
\[
4\, e^{2\la}\vec{H}_{\vec{\xi}}\equiv2\,\Delta\vec{\xi}=\nabla^\perp S\cdot\nabla\vec{\xi}-\nabla\vec{R}\res\nabla^\perp\vec{\xi}\quad.
\]
In particular, since $|\nabla\vec{\xi}|^2=2 e^{\la}$, calling upon (\ref{VII.32}) in the latter yields the announced (\ref{VII.2}), thereby completing the proof.\hfill $\Box$

\section{$L^2$-weak energy quantization for the Gauss map of a Willmore immersion in neck regions}
\reset

\medskip

\begin{Lma}
\label{lm-A.1}{\bf $[L^2-$weak energy quantization for the Gauss map in neck regions]}
There exists a constant $\ep(m)>0$ with the following property. Let $\vec{\xi}_k$ be a sequence of conformal Willmore immersions from $B_{R_k}(0)$ into ${\R}^m$,
with $R_k\rightarrow+\infty$. Whenever $r_k\rightarrow 0$ satisfies
\be
\label{A.1}
\sup_{r\in (r_k,R_k/2)}\int_{B_{2\, r}(0)\setminus B_r(0)}|\nabla\vec{n}_{\vec{\xi}_k}|^2\, dx\le\ep(m)\quad,
\ee
and\footnote{as usual, $\la_k$ is the conformal parameter of $\vec{\xi}_k$.}
\be
\label{A.2}
\|\nabla \la_k\|_{L^{2,\infty}(\Om_k)}+\int_{\Om_k}|\nabla\vec{n}_{\vec{\xi}_k}|^2\, dx\le\La<+\infty\qquad\forall\:k\quad,
\ee
where $\La>0$ is independent of $k$ and $\Om_k=B_{R_k}(0)\setminus B_{r_k}(0)$ ; then
\be
\label{A.4}
\forall\: \ep\in(0,\ep(m))\quad\exists\:\al\in (0,1)\quad\mbox{s.t.}\quad\quad\limsup_{k\rightarrow +\infty}\||x|\,\nabla\vec{n}_{\vec{\xi}_k}(x)\|_{L^\infty(B_{\al R_k/2}(0)\setminus B_{\al^{-1} r_k}(0))}\le\ep\quad.
\ee
Hence, in particular,
\be
\label{A.5}
\forall\: \ep\in(0,\ep(m))\quad\exists\:\al\in (0,1)\quad\mbox{s.t.}\quad\quad\limsup_{k\rightarrow +\infty}\|\nabla\vec{n}_{\vec{\xi}_k}\|_{L^{2,\infty}(B_{\al R_k/2}(0)\setminus B_{\al^{-1} r_k}(0))}\le \ep\quad.
\ee
\hfill$\Box$
\end{Lma}
{\bf Proof of lemma~\ref{lm-A.1}.}
Analogously to the argument given in the paragraph following (\ref{VII.4}) in the proof of lemma~\ref{lm-VII.1}, the constant $\ep(m)$ may be chosen so as to ensure that there holds
\be
\label{A.z6z}
|\nabla\vec{n}_{\vec{\xi}_k}(x)|^2\le C\,|x|^{-2}\int_{B_{2|x|}(0)\setminus B_{|x|/2}(0)}|\nabla\vec{n}_{\vec{\xi}_k}|^2\, dx\le C\,\ep(m)\,x^{-2}\qquad\forall\: x\in B_{R_k/2}(0)\setminus B_{2r_k}(0)\quad,
\ee
where $C$ is some constant independent of the data of the problem. 
This implies in particular that
\be
\label{VIII.00}
\|\nabla\vec{n}_{\vec{\xi}_k}\|_{L^{2,\infty}( B_{R_k/2}(0)\setminus B_{2r_k}(0))}\le C'\,\sqrt{\ep(m)}\quad.
\ee
We will argue by contradiction, by assuming that there exists a sequence of conformal Willmore immersions on $B_{R_k}(0)$ satisfying (\ref{A.1}) and (\ref{A.2}), but for which
there exist $\ep_1>0$ and $x_k\in\Om_k$ with
\be
\label{A.7}
\log\lf|\frac{|x_k|}{r_k}\rg|\rightarrow +\infty\quad\quad\mbox{ , }\quad\quad\log\lf|\frac{|x_k|}{R_k}\rg|\rightarrow -\infty\quad,
\ee
and with
\be
\label{A.8}
|x_k|\,|\nabla\vec{n}_{\vec{\xi}_k}(x_k)|\ge\ep_1>0\quad.
\ee
From (\ref{A.z6z}), we deduce in particular that
\be
\label{A.9}
\int_{B_{2|x_k|}(0)\setminus B_{|x_k|/2}(0)}|\nabla\vec{n}_{\vec{\xi}_k}|^2\, dx\ge\frac{\ep^2_1}{C}>0
\ee
We next demand that $\ep(m)$ be smaller than its counterpart from lemma~\ref{lm-VII.1}. Since all the conditions of this lemma are fulfilled, there exist $\vec{L}_k\in C^\infty(B_{R_k}(0),{\R}^m)$, $S_k\in C^\infty(B_{R_k}(0),{\R})$, and
 $\vec{R}_k\in C^\infty(B_{R_k}(0),\bigwedge^2{\R}^m)$  such that
\be
\label{A.10}
\left\{
\begin{array}{lcl}
\nabla^\perp\vec{L}_k&:=&\nabla\vec{H}_k\,-\,3\,\pi_{\vec{n}_{\vec{\xi}_k}}\big(\nabla\vec{H}_k\big)\,+\,\star\,\big(\nabla^\perp\vec{n}_{\vec{\xi}_k}\wedge\vec{H}_k\big)\\[3mm]
\nabla S_k&:=&\vec{L}_k\cdot\nabla\vec{\xi}_k\\[3mm]
\nabla\vec{R}_k&:=&\vec{L}_k\wedge\nabla\vec{\xi}_k\,+\,2\,\vec{H}_k\wedge\nabla^\perp\vec{\xi}_k\quad.
\end{array}
\right.
\ee

\noindent
In the course of the proof of lemma~\ref{lm-VII.1}, we have seen that
\be
\label{VIII.1}
\| e^{\la_k}\vec{L}_{\vec{\xi}_k}\|_{L^{2,\infty}(B_{R_k/2}\setminus B_{2r_k})}+\|\nabla S_k\|_{L^{2,1}(B_{R_k/2}\setminus B_{2r_k})}+\|\nabla \vec{R}_k\|_{L^{2,1}(B_{R_k/2}\setminus B_{2r_k})}\le C(m,\La)\quad.
\ee
Moreover, we have the conservative conformal Willmore system
\be
\label{VIII.2}
\left\{
\begin{array}{lcl}
\ds\Delta S_k&=&\ds-(\star\nabla\vec{n}_{\vec{\xi}_k})\cdot\nabla^\perp\vec{R}_k\\[3mm]
\ds\Delta\vec{R}_k&=&\ds(-1)^m \star(\nabla\vec{n}_{\vec{\xi}_k}\bullet\nabla^\perp\vec{R}_k)+(\star\nabla\vec{n}_{\vec{\xi}_k})\cdot\nabla^\perp S_k\\[3mm]
2\Delta\vec{\xi}_k&=&\nabla^\perp S_k\cdot\nabla\vec{\xi}_k-\nabla\vec{R}_k\res\nabla^\perp\vec{\xi}_k\quad.
\end{array}
\right.
\ee
Consider the conformal mapping
\[
\ti{\vec{\xi}}_k(y):=e^{-\la_k(x_k)}\, (\vec{\xi}_k(|x_k|\,y)-\vec{\xi}_k(x_k))\quad.
\]
An elementary computation shows that
\[
\nabla_y\vec{n}_{\ti{\vec{\xi}}_k}(y)=|x_k|\ \nabla_x\vec{n}_{\vec{\xi}_k}(|x_k|\, y)\quad,
\]
and
\[
\nabla_y\ti{\vec{H}}_k(y)=|x_k|\ e^{-\la_k(x_k)}\ \nabla_x\vec{H}_k(|x_k|\, y)\quad,
\]
where $\ti{\vec{H}}_k$ is the mean curvature vector of the immersion $\ti{\vec{\xi}}_k$. The corresponding ${\R}^m$-valued map $\ti{\vec{L}}_k$ satisfies
\[
\nabla_y\ti{\vec{L}}_k(y)=|x_k|\ e^{-\la_k(x_k)}\ \nabla_x\vec{L}_k(|x_k|\, y)\quad,
\]
so that
\[
\ti{\vec{L}}_k(y)=e^{-\la_k(x_k)}\ \vec{L}_k(|x_k|\, y)\quad.
\]
Clearly, the conformal factor of the rescaled immersion satisfies $\ti{\la}_k(y)=\la_k(|x_k|\, y)-\la_k(x_k)$. This implies in particular that the corresponding $\ti{S}_k$ and $\ti{\vec{R}}_k$ are given by
\[
\ti{S}_k(y)=S_k(|x_k|\, y)\quad\quad\mbox{ and }\quad\quad\ti{\vec{R}}_k(y)=\vec{R}_k(|x_k|\, y)\quad.
\]
We thus see that the maps $S_k$ and ${\vec{R}_k}$ behave under rescaling like the Gauss map $\vec{n}_{\vec{\xi}}$\,. This remarkable fact was first observed in \cite{BR2}, where it played a decisive role.  From the latter and (\ref{VIII.1}), it follows that
\be
\label{VIII.2b}
\limsup_{k\rightarrow +\infty}\;\|\nabla \ti{S}_k\|_{L^{2,1}(K)}+\|\nabla \ti{\vec{R}}_k\|_{L^{2,1}(K)}\le C(m,\La)\quad,
\ee
for any compact subdomain $K$ of $\C\setminus\{0\}$.

\medskip

The pointwise control of the conformal factor in neck regions provided by
lemma~\ref{lm-Aa0} shows that the sequence of Willmore immersions $\ti{\vec{\xi}_k}$ satisfies
\be
\label{VIII.3}
\limsup_{k\rightarrow +\infty}\;\|\nabla\vec{n}_{\ti{\vec{\xi}}_k}\|_{L^\infty_{loc}({\C}\setminus \{0\})}<+\infty
\ee
and 
\be
\label{VIII.4}
\limsup_{k\rightarrow +\infty}\;\|\log|\nabla{\ti{\vec{\xi}}_k}|\|_{L^\infty_{loc}({\C}\setminus \{0\})}<+\infty
\ee
From the work in \cite{Ri2}, we conclude that for all $l\in\N$, the sequence $\vec{\xi}_k$ converges in $C^l_{loc}({\C}\setminus\{0\})$ to a Willmore immersion $\ti{\vec{\xi}}_\infty$ of ${\C}\setminus\{0\}$. This strong convergence implies that the condition (\ref{A.9}) and the system (\ref{VIII.2}) pass to the limit. Hence,
\be
\label{VIII.5}
\int_{B_{2}(0)\setminus B_{1/2}(0)}|\nabla\vec{n}_{\ti{\vec{\xi}}_\infty}|^2\, dx\ge\frac{\ep^2_1}{C}>0\quad,
\ee
and
\be
\label{VIII.6}
\left\{
\begin{array}{rcl}
\ds\Delta \ti{S}_\infty&=&\ds-(\star\nabla\vec{n}_{\ti{\vec{\xi}}_\infty})\cdot\nabla^\perp\ti{\vec{R}}_\infty\\[3mm]
\ds\Delta\ti{\vec{R}}_\infty&=&\ds(-1)^m \star(\nabla\vec{n}_{\ti{\vec{\xi}}_\infty}\!\bullet\nabla^\perp\ti{\vec{R}}_\infty)+(\star\nabla\vec{n}_{\ti{\vec{\xi}}_\infty})\cdot\nabla^\perp\ti{S}_\infty\\[3mm]
4 e^{2\ti{\la}_\infty}\vec{H}_{\ti{\vec{\xi}}_\infty}&=&2\Delta\ti{\vec{\xi}}_\infty\:=\:\nabla^\perp\ti{S}_\infty\cdot\nabla\ti{\vec{\xi}}_\infty-\nabla\ti{\vec{R}}_\infty\res\nabla^\perp\ti{\vec{\xi}}_\infty\quad.
\end{array}
\right.
\ee
\noindent
On one hand, we obtain from (\ref{VIII.00}) that
\be
\label{VIII.7}
\|\nabla\vec{n}_{\ti{\vec{\xi}}_\infty}\|_{L^{2,\infty}({\C})}\le C\, \sqrt{\ep(m)}\quad.
\ee
On the other hand, it follows from (\ref{VIII.2b}) that
\be
\label{VIII.8}
\|\nabla \ti{S}_\infty\|_{L^{2,1}({\C})}+\|\nabla \ti{\vec{R}}_\infty\|_{L^{2,1}({\C})}\le C(m,\La)<+\infty\quad.
\ee
Applying to the system (\ref{VIII.6}) the Wente inequality from lemma~\ref{weakwente} with $(p,q)=(2,1)$ yields now
\be
\label{VIII.9}
\begin{array}{rl}
\|\nabla \ti{S}_\infty\|_{L^{2,1}({\C})}+\|\nabla \ti{\vec{R}}_\infty\|_{L^{2,1}({\C})}&\le C\ \|\nabla\vec{n}_{\ti{\vec{\xi}}_\infty}\|_{L^{2,\infty}({\C})}\ \lf[\|\nabla \ti{S}_\infty\|_{L^{2,1}({\C})}+\|\nabla \ti{\vec{R}}_\infty\|_{L^{2,1}({\C})}\rg]\quad.\\[5mm]
 &\le C_m\, \sqrt{\ep(m)}\ \lf[\|\nabla \ti{S}_\infty\|_{L^{2,1}({\C})}+\|\nabla \ti{\vec{R}}_\infty\|_{L^{2,1}({\C})}\rg]\quad,
\end{array}
\ee
where we have used (\ref{VIII.7}). Accordingly, when $C_m\sqrt{\ep(m)}<1$, we deduce that $\nabla \ti{S}_\infty=0$ and $\nabla\ti{\vec{R}}_\infty=0$. Injecting this information in the last identity of (\ref{VIII.6}) gives
\[
\vec{H}_{\ti{\vec{\xi}}_\infty}\equiv \vec{0}\quad.
\]
In particular, since the (weighted) mean curvature vector is the trace part of the gradient of the Gauss map, it follows easily that
\be
\label{VIII.9a}
|\nabla \vec{n}_{\ti{\vec{\xi}}_\infty}|^2=-2\,e^{2\ti{\la}_\infty}K_{\ti{\vec{\xi}}_\infty}\, \quad\quad\mbox{ on }{\C}\quad,
\ee
where $K_{\ti{\vec{\xi}}_\infty}$ is the Gauss curvature of the limit immersion ${\ti{\vec{\xi}}_\infty}$.\\
Let us next choose the constant $\ep(m)$ small enough that we can apply lemma~\ref{lm-AA.1} with the condition (\ref{VIII.7}). Namely, on any ball $B_\rho(0)$, we construct a moving frame $\{\bbe^{\,\rho}_1,\bbe^{\,\rho}_2\}$ satisfying
 \[
\star\,\hat{\bn}_{\ti{\vec{\xi}}_\infty}\,=\,{\bbe^{\,\rho}_1}\wedge{\bbe^{\,\rho}_2}\quad,\qquad div\,({\bbe^{\,\rho}_1}\cdot\nabla{\bbe^{\,\rho}_2})\,=\,0
\]
and 
\be
\label{VIII.10}
\int_{B_R(0)}|\nabla\bbe^{\,\rho}_1|^2+|\nabla\bbe^{\,\rho}_2|^2\,\leq\,C\,\int_{{\C}}|\nabla{\bn}_{\ti{\vec{\xi}}_\infty}|^2\quad.
\ee
We now consider a sequence of radii $\rho_j$ converging to infinity. Extracting a subsequence gives the existenceof a limiting frame $(\bbe^{\,\infty}_1,\bbe^{\,\infty}_2)$ on ${\C}$ satisfying
 \[
\star\,\hat{\bn}_{\ti{\vec{\xi}}_\infty}\,=\,{\bbe^{\,\infty}_1}\wedge{\bbe^{\,\infty}_2}\quad,\qquad div({\bbe^{\,\infty}_1}\cdot\nabla{\bbe^{\,\infty}_2})\,=\,0\quad,
\]
and 
\be
\label{VIII.11}
\int_{{\C}}|\nabla\bbe^{\,\infty}_1|^2+|\nabla\bbe^{\,\infty}_2|^2\,\leq\,C\,\int_{{\C}}|\nabla{\bn}_{\ti{\vec{\xi}}_\infty}|^2\quad.
\ee
We have at every point of the {\it minimal} conformal immersion $\ti{\vec{\xi}}_\infty$ the identity
\[
 e^{2\ti{\la}_\infty}K_{\ti{\vec{\xi}}_\infty}=-\,\nabla^\perp\bbe^{\,\infty}_1\cdot\nabla\bbe^{\,\infty}_2\quad.
\]
Integrating over ${\C}$ then yields
\[
\int_{{\C}}K_{\ti{\vec{\xi}}_\infty}\ e^{2\ti{\la}_\infty}\ dx=\int_{{\C}} d\bbe^{\,\infty}_1\dot{\wedge} d\bbe^{\,\infty}_2=0
\]
Paired to (\ref{VIII.9a}), the latter implies that $\nabla\vec{n}_{\ti{\vec{\xi}}_\infty}\equiv 0$, which contradicts (\ref{VIII.5}). This is the desired contradiction, and thus lemma~\ref{lm-A.1} is proved.\hfill $\Box$

\section{Proof of theorems}
\reset
\subsection{Proof of theorem~\ref{th-I.1}}

Let $\vec{\Phi}_k$ be a sequence of Willmore immersions of $\Sigma$ with uniformly bounded energy 
\[
\limsup_{k\rightarrow +\infty}\; E(\vec{\Phi}_k)<+\infty\quad,
\]
and such that the conformal class of the induced metric $g_k=\vec{\Phi}_k^\ast g_{{\R}^m}$ remains in a compact subdomain of the moduli space of $\Sigma$.

\medskip

We consider a subsequence, still denoted $\vec{\Phi}_k$, as given by lemma~\ref{lm-I.1}. Let $\Xi_k$ be the corresponding sequence of M\"obius transformations, and let $f_k$ be the sequence of Lipschitz diffeomorphisms of $\Sigma$, as given by lemma~\ref{lm-I.1}. 
We set $\vec{\xi}_k:=\Xi_k\circ\vec{\Phi}_k\circ f_k$. The conclusion of lemma~\ref{lm-I.1} states that $\vec{\xi}_k$ weakly
converges in $W^{2,2}$ to a weak immersion $\vec{\xi}_\infty$, away from finitely many points $a_i$. This limiting immersion might be branched at the points $a_i$.
 
\medskip

Let $h_k$ be constant scalar curvature metric with respect to which $\vec{\xi}_k$ is conformal, and let $u_k$ be the function satisfying
\[
g_k=\vec{\xi}_k^\ast g_{{\R}^m}=e^{2\,u_k}\,h_k\quad.
\]
The classical Liouville equation reads
\[
\Delta_{h_k} u_k+ K_{g_k} e^{2\, u_k}=K_{h_k}\quad.
\]
Since $h_k$ strongly converges to some limiting constant scalar curvature $h_\infty$, classical results from geometric analysis on manifolds (see for instance \cite{Aub}) show that the gradient of the Green function associated with the Laplace-Beltrami operator $\Delta_{h_k}$ is uniformly bounded in $L^{2,\infty}$. Hence,
\be
\label{IX.0}
\|\nabla u_k\|_{L^{2,\infty}(\Sigma)}\le C\, \lf[\int_{\Sigma}|K_{g_k}|\ dvol_{g_k}+\int_{\Sigma}|K_{h_k}|\ dvol_{h_k}\rg]\le C\, \lf[E(\vec{\Phi}_k)+2\pi\,|\chi(\Sigma)|\rg]\quad,
\ee
where the $L^{2,\infty}$-norm is to be understood with respect to a reference metric $g_0$ on $\Sigma$. As before, $\chi(\Sigma)$ denotes the Euler characteristic of $\Sigma$.
We define the constant
\be
\label{IX.1a}
\La:=\sup_{k\in{\N}} E(\vec{\Phi}_k)+\|\nabla u_k\|_{L^{2,\infty}(\Sigma)}<+\infty\quad.
\ee

\medskip

Owing to the $\ep$-regularity theorem I.5 in \cite{Ri2}, the sequence $\vec{\xi}_k$ actually converges strongly in $C^l_{loc}(\Sigma\setminus\{a_1,\ldots,  a_N\},{\R}^m)$,
for all $l\in{\N}$. Namely,
\be
\label{IX.0a}
\vec{\xi}_k\longrightarrow\vec{\xi}_\infty\quad\quad\mbox{ in }\quad\quad C^l_{loc}(\Sigma\setminus\{a_1,\ldots,  a_N\},{\R}^m)\quad.
\ee
Accordingly, the limiting immersion $\vec{\xi}_\infty$ is Willmore, although a priori only away from the points $\{a_1,\ldots, a_N\}$. Since $\vec{\xi}_k$
is Willmore, for any fixed small radius $\rho>0$ and any $i=1,\ldots, N$, one has
\[
\int_{B_\rho(a_i)}div\lf( 2\,\nabla\vec{H}_{\vec{\xi}_k}-3\,\nabla(\pi_{\vec{n}_{\vec{\xi}_k}})(\vec{H}_{\vec{\xi}_k})-\star(\nabla^\perp\vec{n}_{\vec{\xi}_k}\wedge\vec{H}_{\vec{\xi}_k})\rg)\ dx=0\quad,
\]
where $B_\rho(a_i)$ is the geodesic ball centered at $a_i$ and of radius $\rho$ for the flat metric, in some converging conformal coordinate system for $\vec{\xi}_k^\ast g_{{\R}^m}$. Whence,
\[
\int_{\p B_\rho(a_i)} 2\,\p_\nu\vec{H}_{\vec{\xi}_k}-3\,\p_{\nu}(\pi_{\vec{n}_{\vec{\xi}_k}})(\vec{H}_{\vec{\xi}_k})+\star(\p_\tau\vec{n}_{\vec{\xi}_k}\wedge\vec{H}_{\vec{\xi}_k}))\ \rho\, d\theta=0\quad,
\]
where $\nu$ is the unit normal vector $x/|x|$ and $\tau$ the unit tangent vector $x^\perp/|x|$ to the boundary of $B_\rho(a_i)$. The strong convergence of $\vec{\xi}$ to $\vec{\xi}_{\infty}$ enables passing to the limit in the last identity, and we thus find
\[
\int_{\p B_\rho(a_i)} 2\,\p_\nu\vec{H}_{\vec{\xi}_\infty}-3\,\p_{\nu}(\pi_{\vec{n}_{\vec{\xi}_\infty}})(\vec{H}_{\vec{\xi}_\infty})+\star(\p_\tau\vec{n}_{\vec{\xi}_\infty}\wedge\vec{H}_{\vec{\xi}_\infty}))\ \rho\, d\theta=0\quad,
\]
for all fixed $\rho$. \\
This is precisely the condition required to apply the ``point removability" theorem established in \cite{Ri2}. In particular, $\vec{\xi}_\infty$ extends to a (possibly branched), smooth Willmore immersion in $B_\rho(a_i)$. More detailed information on the behavior of the immersion near the points $a_i$ are given in \cite{BR2}.

\medskip

We now choose a fixed $\ep_0>0$ smaller than the constants $\ep(m)$ of lemmas~\ref{lm-VII.1} and ~\ref{lm-A.1}. For this particular choice of $\ep_0$, we apply the bubble-neck decomposition procedure outlined in proposition~\ref{pr-III.1}. Owing to (\ref{III.12}), both lemma~\ref{lm-VII.1} and lemma~\ref{lm-A.1} apply in each connected component
of $\Om_k(\al)$, which, recall, is a union of disjoint annuli. From lemma~\ref{lm-VII.1}, we thus have
\be
\label{IX.1}
\limsup_{k\rightarrow +\infty}\|e^\la\vec{H}_{\vec{\xi}_k}\|_{L^{2,1}(\Om_k(1/2))}\le C(m,\La)\quad.
\ee
As for lemma~\ref{lm-A.1}, it yields
\be
\label{IX.2}
\forall\:\ep>0\quad\quad\exists\,\al>0\quad\quad\mbox{ s.t. }\quad\quad\|\nabla\vec{n}_{\vec{\xi}_k}\|_{L^{2,\infty}(\Om_k(\al))}\le \ep\quad.
\ee
Combining this last assertion with the result of lemma~\ref{lm-A.a00} implies that
\be
\label{IX.3}
\forall\:\ep>0\quad\quad\exists\,\al>0\quad\quad\mbox{ s.t. }\quad\quad\lf|\int_{\Om_k(\al)} K_{\vec{\xi}_k}\ dvol_{g_k}\rg|\le\ep\quad.
\ee
Altogether now, (\ref{IX.1}) and (\ref{IX.2}) give
\be
\label{IX.4}
\begin{array}{l}
\ds\forall\:\ep>0\quad\quad\exists\,\al>0\quad\exists\, k_0\in{\N}\quad\mbox{ s.t. }\quad\forall\: k\ge k_0\\[5mm]
\ds\|e^\la\vec{H}_{\vec{\xi}_k}\|_{L^1(\Om_k(\al))}\le \|e^\la\vec{H}_{\vec{\xi}_k}\|^{1/2}_{L^{2,1}(\Om_k(\al))}\ \|e^\la\vec{H}_{\vec{\xi}_k}\|^{1/2}_{L^{2,\infty}(\Om_k(\al))}\le \sqrt{C(m,\La)\ \ep}\quad,
\end{array}
\ee
where we used the fact that $L^{2,\infty}$ is the dual space of $L^{2,1}$. \\
Combining (\ref{IX.4}) and (\ref{IX.3}) shows that
\be
\label{IX.5}
\lim_{\al\rightarrow 0}\,\lim_{k\rightarrow +\infty}\int_{\Om_k(\al)}|\nabla\vec{n}_{\vec{\xi}_k}|^2\, dx=0\:.
\ee
This is the {\it no-neck energy property} which we aimed at establishing. To finish the proof of theorem~\ref{th-I.1}, there remains to identify the energy of the bubbles as converging in the limit to that of Willmore spheres. This is what we do below. 

\medskip

Let $i\in\{1,\ldots, N\}$ and let $j\in\{1,\ldots, Q^i\}$. From (\ref{III.14}), on the bubble domain
\[
B(i,j,\al,k):=B_{\al^{-1}\rho_k^{i,j}}(x^{i,j}_k)\setminus \cup_{j'\in I^{i,j}}\,B_{\al\rho_k^{i,j}}(x^{i,j'}_k)\quad,
\]
there is no concentration of energy. Hence, arguing as in the proof of lemma~\ref{lm-A.0} (see also lemma III.1 in \cite{Ri3}), the conformal factor satisfies a uniform Harnack estimate:
\be
\label{IX.6}
\forall\: \, 0<\al<1\quad\exists\; C_\al>0\quad\quad\mbox{ s.t. }\quad\quad1\le\frac{\sup_{\{x\in B(i,j,\al,k)\}} e^{\la_k}(x)}{\inf_{\{x\in B(i,j,\al,k)\}} e^{\la_k}(x)}\le C_\al\quad.
\ee
Choose an arbitrary point 
\[
z^{i,j}_k\in B(i,j,\al,k)\cap B_{\rho_k^{i,j}}(x_k^{i,j})\setminus \lf(\cup_{j'\in I^{i,j}}\,B_{M^{-1}\rho_k^{i,j}}(x^{i,j'}_k)\rg)\quad,
\]
where $M>$card$(I^{i,j})$\,; and set $\la(i,j,k,\al):= {\la_k}(x_k)$. We introduce the renormalized Willmore immersion 
\[
\ti{\vec{\xi}}_k(y):= e^{-\la(i,j,\al,k)}\ \lf[\vec{\xi}_k(\rho_k^{i,j}\ y+x_k^{i,j})-\vec{\xi}_k(z^{i,j}_k)\rg]\quad.
\]
We can extract a subsequence in such a way that the following limits exist:
\[
\forall\:\ j'\in I^{i,j}\quad\quad \lim_{k\rightarrow +\infty} \frac{x_k^{i,j}-z_k^{i,j}}{\rho^{i,j}_k}=a^{j,j'}_i\in B_1(0)\setminus\{0\}\quad.
\]
As the energy does not concentrate on any $B(i,j,\al,k)$ for $0<\al<1$, we obtain, using once more the $\ep$-regularity and the local control of the conformal factor ensured by (\ref{IX.6}), that 
\[
\ti{\vec{\xi}}_k\longrightarrow\ti{\vec{\xi}}_\infty\quad\quad\mbox{ in }\quad C^l_{loc}({\C}\setminus\cup_{j'\in I^{i,j}}\{a^{j,j'}_i\})\quad\quad\forall\:\ l\in{\N}\quad,
\]
where $\ti{\vec{\xi}}_\infty$ is a Willmore immersion of ${\C}\setminus\cup_{j'\in I^{i,j}}\{a^{j,j'}_i\}$. Suppose that
\[
\int_{{\C}}|\nabla\ti{\vec{\xi}}_\infty|^2\, dx<+\infty\quad.
\]
Then, using the stereographic projection of $S^2$ onto ${\C}$ with respect to the north pole, $\ti{\vec{\xi}}_\infty$ realizes an immersion of $S^2\setminus\cup_{j'\in I^{i,j}}\{a^{j,j'}_i\}\cup\{North\}$ included in a bounded domain. Just as we explained above for the limit immersion $\vec{\xi}_\infty$, one can remove the singularities $\cup_{j'\in I^{i,j}}\{a^{j,j'}_i\}\cup\{North\}$, which possibly realize branch points of our Willmore $S^2$.
We have moreover
\be
\label{IX.6a}
\lim_{\al\rightarrow 0}\lim_{k\rightarrow +\infty}\int_{B(i,j,\al,k)} e^{2\la_k}\ |\vec{H}_{\vec{\xi}_k}|^2=W(\ti{\vec{\xi}}_\infty)\quad.
\ee
On the other hand, suppose that\footnote{This should correspond to the case when $\la(i,j,\al,k)\rightarrow -\infty$.} 
\[
\int_{{\C}}|\nabla\ti{\vec{\xi}}_\infty|^2\ dx_1\, dx_2=+\infty\quad.
\]
Following \cite{Ri4}, we can find a point $p_k\in{\R}$ such that
\[
\forall\:\, k\in {\N}\quad\quad e^{-\la(i,j,\al,k)}\lf[\vec{\xi}_k(\Sigma)-\vec{\xi}_k(z^{i,j}_k)\rg]\cap B_1(p_k)=\emptyset\quad\quad\mbox{ and }\quad\quad\limsup_{k\rightarrow+\infty}|p_k|<+\infty\quad.
\]
Let ${\mathcal I}_k$ be the inversion with respect to $p_k$, and $\hat{\vec{\xi}_k}:={\mathcal I}_k\circ\ti{\vec{\xi}}_k$. We may assume, modulo extraction of a subsequence if necessary, that the sequence $p_k$ converges to a point $p_\infty$. We let ${\mathcal I}_\infty$ be the inversion in ${\R}^m$ with respect to this limit-point. There holds
\[
\hat{\vec{\xi}}_k\longrightarrow\hat{\vec{\xi}}_\infty={\mathcal I}_\infty\circ\ti{\vec{\xi}}_\infty\quad\quad\mbox{ in }\quad C^l_{loc}({\C}\setminus\cup_{j'\in I^{i,j}}\{a^{j,j'}_i\})\quad\quad\forall\:\ l\in{\N}\quad,
\]
where $\hat{\vec{\xi}}_\infty$ is a Willmore immersion of ${\C}\setminus\cup_{j'\in I^{i,j}}\{a^{j,j'}_i\}$ satisfying
\[
\int_{{\C}}|\nabla\hat{\vec{\xi}}_\infty|^2\, dx<+\infty\quad.
\]
Arguing as above, $\hat{\vec{\xi}}_\infty$ extends to a smooth, possibly branched, immersion of $S^2$. Using lemma A.4 in \cite{Ri4}, we then obtain that
\be
\label{IX.7}
\lim_{\al\rightarrow 0}\lim_{k\rightarrow +\infty}\int_{B(i,j,\al,k)} e^{2\la_k}\ |\vec{H}_{\vec{\xi}_k}|^2\, dx=W(\hat{\vec{\xi}}_\infty)-4\pi\,\theta_0\quad,
\ee
where $\theta_0$ is the integer density of ${\vec{\xi}}_\infty({\C})$ at the image point $0\in{\R}^m$. 

\medskip

Altogether, (\ref{IX.0a}), (\ref{IX.6}), and (\ref{IX.7}), imply that there exists a subsequence (indexed as the original, for notational convenience) such that
\be
\label{IX.8}
\lim_{k\rightarrow +\infty}W(\vec{\Phi}_k)=W(\vec{\xi}_\infty)+\sum_{s=1}^p W(\ti{\vec{\xi}^s}_\infty)+\sum_{t=1}^q \lf[W(\hat{\vec{\xi}^t}_\infty)-4\pi\theta_0^t\rg]\quad.
\ee
The index $s$ is used for the first alternative, in which the bubble has no infinite end ; while the index $t$ is used for the second alternative, in which the bubble has at least one infinite
planar end, and it must thus be inverted to become compact. This is exactly the desired {\it energy identity} (\ref{0.2}), thereby ending the proof of theorem~\ref{th-I.1}.\hfill $\Box$

\subsection{Proof of theorem~\ref{th-I.3}}

Let $\vec{\Phi}_k$ be a sequence of Willmore immersions satisfying
$$
W(\vec{\Phi}_k)<\min\{8\pi,\om^m_g\}-\delta\quad.
$$

We want to show that it is possible extract a subsequence $\vec{\Phi}_{k'}$, and to find a sequence of diffeomorphisms $f_{k'}$ of $\Sigma$ and Moebius transformationa $\Xi_{k'}$, such that
\be
\label{X.1}
\vec{\xi}_{k'}:=\Xi_{k'}\circ\vec{\Phi}_{k'}\circ f_{k'}\longrightarrow \vec{\xi}_\infty\quad\quad \mbox{strongly in}\: C^l(\Sigma)\quad\forall\:l\in\N\quad.
\ee

From the work in \cite{Ri4} (cf. also \cite{KuLi}), the sequence of conformal classes of the metrics induced by $\vec{\Phi}_k$ remains in a compact subdomain of the Moduli space. Theorems \ref{th-I.2} and \ref{pr-III.4} whence apply.\\
If
\[
\lim_{k'\rightarrow\infty} W(\vec{\Phi}_{k'})=W(\vec{\xi}_\infty)\quad.
\]
then from theorem~\ref{th-I.2}, we get (\ref{X.1}) directly. If not, then there must be a bubble. We take it to be the most concentrated one, centered at $x_{k'}$ and of radius $\rho_{k'}$, and converging to a point $a^1$.
As is explained at the end of the previous section, we can find a Moebius transformation $\La_{k'}$ such that in local conformal coordinates there holds
\[
\La_{k'}\circ\vec{\xi}_{k'}(\rho_{k'}y+x_{k'})\longrightarrow \hat{\vec{\xi}}_{\infty}\quad,
\]
where $\hat{\vec{\xi}}_{\infty}$ is a Willmore sphere. Assume this sphere is homothetic to the round sphere.
Then, if $\La_{k'}$ contains no inversion and $\hat{\vec{\xi}}_\infty$ is of the type $\eta_s$ from theorem~\ref{th-I.1}, we conclude that
\[
\lim_{k'\rightarrow\infty} W(\vec{\Phi}_{k'})\ge W(\vec{\xi}_\infty) +W(\hat{\vec{\xi}}_{\infty})\ge 4\pi +4\pi=8\pi\quad,
\]
which is a contradiction.\\
If now $\La_{k'}$ contains an inversion, then modulo a dilation in the image, we can arrange for $\xi_{k'}(\rho_{k'}y+x_{k'})$ to converge strongly to a plane, and thus
\[
\int_{B_{\rho_{k'}}(x_{k'})}|\nabla{\vec{\xi}_{k'}}|\rightarrow 0\quad.
\]
This contradicts the very definition of a bubble, which requires that
\[
\int_{B_{\rho_{k'}}(x_{k'})}|\nabla{\vec{\xi}_{k'}}|\ge \ep(m)>0\quad.
\]
Accordingly, $\hat{\vec{\xi}}_{\infty}$ can not be homothetic to the round sphere. Calling upon the results from \cite{Bry} and \cite{Mon}, it follows that
\[
\lim_{k'\rightarrow\infty} W(\vec{\Phi}_{k'})= \lim_{k'\rightarrow\infty} W(\La_{k'}\circ\vec{\xi}_{k'})\ge\lim_{k'\rightarrow\infty} W(\La_{k'}\circ\xi_{k'}(\rho_{k'}y+x_{k'}))\longrightarrow W(\hat{\vec{\xi}}_{\infty})\ge 8\pi\:,
\]
which is again a contradiction.

\medskip

Therefore, we have shown that in all cases, the assertion (\ref{X.1}) holds, thereby proving theorem~\ref{th-I.3}.

\end{document}